\documentclass[11pt,english,vecarrow]{extarticle}
\usepackage[T1]{fontenc}
\usepackage[utf8]{inputenc}
\usepackage[letterpaper]{geometry}
\usepackage{babel}
\usepackage{array}
\usepackage{verbatim}
\usepackage{float}
\usepackage{calc}
\usepackage{framed}
\usepackage{amsmath}
\usepackage{amssymb}
\usepackage{stmaryrd}
\usepackage{graphicx}
\usepackage{setspace}
\usepackage{wasysym}

\makeatletter

\providecommand{\tabularnewline}{\\}
\newcommand{\lyxdot}{.}

\numberwithin{equation}{section}
\numberwithin{figure}{section}
\numberwithin{table}{section}

\usepackage{bbm}
\usepackage{tensind}
\tensordelimiter{?}
\usepackage{multicol}
\usepackage[labelfont={sc,small},textfont=small]{caption}

\makeatother

\begin{document}
\title{Gauge Theories and Fiber Bundles:\\
Definitions, Pictures, and Results\\
\quad{}}
\author{\doublespacing{}\textsc{Adam Marsh}}
\date{\textit{\normalsize{}August 17, 2022}}
\maketitle
\begin{abstract}
A pedagogical but concise overview of fiber bundles and their connections
is provided, in the context of gauge theories in physics. The emphasis
is on defining and visualizing concepts and relationships between
them, as well as listing common confusions, alternative notations
and jargon, and relevant facts and theorems. Special attention is
given to detailed figures and geometric viewpoints, some of which
would seem to be novel to the literature. Topics are avoided which
are well covered in textbooks, such as historical motivations, proofs
and derivations, and tools for practical calculations. The present
paper is best read in conjunction with the similar paper on Riemannian
geometry cited herein. 
\end{abstract}
\tableofcontents{}

\section{Introduction}

A manifold includes a tangent space associated with each point. A
frame defines a basis for the tangent space at each point, and a connection
allows us to compare vectors at different points, leading to concepts
including the covariant derivative and curvature. All of these concepts,
covered in a similar style in \cite{Marsh}, can be applied to an
arbitrary vector space associated with each point in place of the
tangent space. This is the idea behind gauge theories. Both manifolds
with connection and gauge theories can then be described using the
mathematical language of fiber bundles.

Throughout the paper, warnings concerning a common confusion or easily
misunderstood concept are separated from the core material by boxes,
as are intuitive interpretations or heuristic views that help in understanding
a particular concept. Quantities are written in \textbf{bold} when
first mentioned or defined.

\section{Gauge theory}

\subsection{\label{subsec:Matter-fields-and-gauges}Matter fields and gauges}

Gauge theories associate each point $x$ on the spacetime manifold
$M$ with a (usually complex) vector space $V_{x}\cong\mathbb{C}^{n}$,
called the \textbf{internal space}\index{internal space}. A $V$-valued
0-form $\vec{\Phi}$ on $M$ is called a \textbf{matter field} (AKA
particle field\index{particle field}). A matter field lets us define
analogs of the quantities associated with a change of frame (see \cite{Marsh})
as follows.

A basis for each $V_{x}$ is called a \textbf{gauge}\index{gauge},
and is the analog of a frame; choosing a gauge is sometimes called
\textbf{gauge fixing}\index{gauge fixing}. Like the frame, a gauge
is generally considered on a region $U\subseteq M$. The analog of
a change of frame is then a (local) \textbf{gauge transformation}\index{gauge transformation}\index{local gauge transformation}
(AKA gauge transformation of the second kind\index{gauge transformation of the second kind}),
a change of basis for each $V_{x}$ at each point $x\in U$. This
is viewed as a representation of a \textbf{gauge group}\index{gauge group}
(AKA symmetry group\index{symmetry group}, structure group\index{structure group})
$G$ acting on $V$ at each point $x\in U$, so that we have 
\begin{equation}
\begin{aligned}\gamma^{-1}\colon & U\to G\\
\rho\colon & G\to GL(V)\\
\Rightarrow\check{\gamma}^{-1}\equiv\rho\gamma^{-1}\colon & U\to GL(V),
\end{aligned}
\end{equation}
and if we choose a gauge it can thus be associated with a matrix-valued
0-form or tensor field 
\begin{equation}
(\gamma^{-1})^{\beta}{}_{\alpha}\colon U\to GL(n,\mathbb{C}),
\end{equation}
so that the components of the matter field $\Phi^{\alpha}$ transform
according to 
\begin{equation}
\Phi^{\prime\beta}=\gamma^{\beta}{}_{\alpha}\Phi^{\alpha}.
\end{equation}
Since all reps of a compact $G$ are similar to a unitary rep, for
compact $G$ we can then choose a \textbf{unitary gauge}\index{unitary gauge},
which is defined to make gauge transformations unitary, so 
\begin{equation}
\check{\gamma}^{-1}\colon U\to U(n).
\end{equation}
This is the analog of choosing an orthonormal frame, where a change
of orthonormal frame then consists of a rotation at each point. A
\textbf{global gauge transformation}\index{global gauge transformation}
(AKA gauge transformation of the first kind\index{gauge transformation of the first kind})
is a gauge transformation that is the same at every point. If the
gauge group is non-abelian (i.e. most groups considered beyond $U(1)$),
the matter field is called a \textbf{Yang-Mills field}\index{Yang-Mills field}
(AKA YM field\index{YM field}). 

\noindent %
\begin{framed}%
\noindent $\triangle$ The term ``gauge group'' can refer to the
abstract group $G$, the matrix rep of this group within $GL(V)$,
the matrix rep within $U(n)$ under a unitary gauge, or the infinite-dimensional
group of maps $\gamma^{-1}$ under composition. This last is sometimes
called the \textbf{global gauge group}\index{global gauge group},
with $G$ or its reps then called the \textbf{local gauge group}\index{local gauge group}.\end{framed}

\noindent %
\begin{framed}%
\noindent $\triangle$ As with vector fields, the matter field $\vec{\Phi}$
is considered to be an intrinsic object, with only the components
$\Phi^{\alpha}$ changing under gauge transformations. \end{framed}

\noindent %
\begin{framed}%
\noindent $\triangle$ Unlike with the frame, whose global existence
is determined by the topology of $M$, there can be a choice as to
whether a global gauge exists or not. This is the essence of fiber
bundles, as we will see in Section \ref{subsec:Fiber-bundles}. \end{framed}

\subsection{The gauge potential and field strength}

We can then define the parallel transporter for matter fields to be
a linear map 
\begin{equation}
\parallel_{C}\colon V_{p}\to V_{q},
\end{equation}
where $C$ is a curve in $M$ from $p$ to $q$. Choosing a gauge,
the parallel transporter can be viewed as a (gauge-dependent) map
\begin{equation}
\parallel^{\beta}{}_{\alpha}\colon\left\{ C\right\} \to GL(n,\mathbb{C}).
\end{equation}
This determines the (gauge-dependent) matter field connection 1-form
\begin{equation}
\Gamma^{\beta}{}_{\alpha}\left(v\right)\colon T_{x}M\to gl(n,\mathbb{C}),
\end{equation}
which can also be written when acting on a $\mathbb{C}^{n}$-valued
0-form as $\check{\Gamma}\left(v\right)\vec{\Phi}$. The values of
the parallel transporter are again viewed as a rep of the gauge group
$G$, so that the values of the connection are a rep of the Lie algebra
$\mathfrak{g}$, and if $G$ is compact we can choose a unitary gauge
so that $\mathfrak{g}$ is represented by anti-hermitian matrices.
We then define the \textbf{gauge potential}\index{gauge potential}
(AKA gauge field\index{gauge field}, four-potential\index{four-potential},
four-vector potential\index{four-vector potential}, vector potential\index{vector potential})
$\check{A}$ by 
\begin{equation}
\check{\Gamma}\equiv-iq\check{A},
\end{equation}
where $q$ is called the \textbf{coupling constant}\index{coupling constant}
(AKA charge\index{charge}, interaction constant\index{interaction constant},
gauge coupling parameter\index{gauge coupling parameter}). Note that
$A^{\beta}{}_{\alpha}$ are then hermitian matrices in a unitary gauge.
The covariant derivative is then 
\begin{equation}
\mathrm{D}_{v}\Phi=\mathrm{d}\vec{\Phi}\left(v\right)-iq\check{A}\left(v\right)\vec{\Phi},
\end{equation}
which is sometimes called the \textbf{gauge covariant derivative}\index{gauge covariant derivative},
and can be generalized to $\mathbb{C}^{n}$-valued $k$-forms as an
exterior covariant derivative 
\begin{equation}
\mathrm{D}\vec{\Phi}=\mathrm{d}\vec{\Phi}-iq\check{A}\wedge\vec{\Phi}.
\end{equation}
For a matter field (0-form), this is often written after being applied
to $e_{\mu}$ as 
\begin{equation}
\mathrm{D}_{\mu}\vec{\Phi}=\partial_{\mu}\vec{\Phi}-iq\check{A}{}_{\mu}\vec{\Phi},
\end{equation}
where $\mu$ is then a spacetime index and 
\begin{equation}
\check{A}_{\mu}\equiv\check{A}(e_{\mu})
\end{equation}
are $gl(n,\mathbb{C})$-valued components. This expression is not
coordinate-dependent, and we may therefore also write it using an
abstract index.

This connection defines a curvature 
\begin{equation}
\check{R}\equiv\mathrm{d}\check{\Gamma}+\check{\Gamma}\wedge\check{\Gamma},
\end{equation}
which lets us define the \textbf{field strength}\index{field strength}
(AKA gauge field\index{gauge field}) $\check{F}$ by 
\begin{equation}
\begin{aligned}\check{R} & \equiv-iq\check{F}\\
\Rightarrow\check{F} & =d\check{A}-iq\check{A}\land\check{A}.
\end{aligned}
\end{equation}
\begin{framed}%
\noindent $\triangle$ The definition $\check{\Gamma}\equiv-iq\check{A}$
is the convention with a mostly pluses metric signature; with a mostly
minuses signature the sign is reversed. However, one also finds this
definition in terms of an elementary charge $e\equiv\pm q$, which
may be positive or negative depending on convention, again reversing
the sign.\end{framed}

\subsection{Spinor fields}

A matter field can also transform as a spinor, in which case it is
called a \textbf{spinor matter field}\index{spinor matter field}
(AKA spinor field\index{spinor field}), and is a 0-form on $M$ which
e.g. for Dirac spinors takes values in $V\otimes\mathbb{C}^{4}$.
The gauge component then responds to gauge transformations, while
the spinor component responds to changes of frame. Similarly, a matter
field on $M^{r+s}$ taking values in $V\otimes\mathbb{R}^{r+s}$ is
called a \textbf{vector matter field}\index{vector matter field}
(AKA vector field\index{vector field}), where the vector component
responds to changes of frame. Finally, a matter field without any
frame-dependent component is called a \textbf{scalar matter field}\index{scalar matter field}
(AKA scalar field\index{scalar field}), and a matter field taking
values in $\mathbb{C}$ (which can be viewed as either vectors or
scalars) is called a \textbf{complex scalar matter field}\index{complex scalar field}
(AKA complex scalar field\index{complex scalar field}, scalar field).
A spinor matter field with gauge group $U(1)$ is called a \textbf{charged
spinor field}\index{charged spinor field}. %
\begin{framed}%
\noindent $\triangle$ It is important remember that spinor and vector
matter fields use the tensor product, not the direct sum, and therefore
cannot be treated as two independent fields. In particular, the field
value $\phi\otimes\psi\in V\otimes\mathbb{C}^{4}$ is identical to
the value $-\phi\otimes-\psi$, which has consequences regarding the
existence of global spinor fields, as we will see in Section \ref{subsec:Spinor-bundles}.\end{framed}

In order to directly map changes of frame to spinor field transformations,
one must use an orthonormal frame so all changes of frame are rotations.
The connection associated with an orthonormal frame is therefore called
a \textbf{spin connection}\index{spin connection}, and takes values
in $so(3,1)\cong\mathrm{spin}(3,1)$. Thus the spin connection and
gauge potential together provide the overall transformation of a spinor
field under parallel transport. All of the above can be generalized
to arbitrary dimension and signature.

\begin{table}[H]
\begin{tabular*}{1\columnwidth}{@{\extracolsep{\fill}}|>{\raggedright}p{0.27\columnwidth}|>{\raggedright}p{0.4\columnwidth}|>{\raggedright}p{0.24\columnwidth}|}
\hline 
Tangent space $T_{p}M=\mathbb{R}^{(r+s)}$ & Spinor space

$S_{p}=\mathbb{K}^{m}$ & Internal space

$V_{p}=\mathbb{C}^{n}$\tabularnewline
\hline 
\hline 
Frame & Standard basis of $\mathbb{K}^{m}$ identified with an initial orthonormal
frame on $M$ & Gauge\tabularnewline
\hline 
Change of frame & $U_{p}\in\mathrm{Spin}(r,s)$ associated with change of orthonormal
frame $\check{\gamma}_{p}$ & Gauge 

transformation\tabularnewline
\hline 
Vector field

$p\mapsto w\in T_{p}M$ & Spinor field

$p\mapsto\psi\in S_{p}$ & Complex / YM field

$p\mapsto\phi\in V_{p}$\tabularnewline
\hline 
Connection

$v\mapsto\check{\Gamma}\left(v\right)\in gl(r,s)$ & Spin connection

$v\mapsto\check{\omega}\left(v\right)\in so(r,s)$, the bivectors & Gauge potential

$v\mapsto\check{A}\left(v\right)\in gl(n,\mathbb{C})$\tabularnewline
\hline 
Curvature

$\check{R}=\mathrm{d}\check{\Gamma}+\check{\Gamma}\wedge\check{\Gamma}$ & Curvature

$\check{R}=\mathrm{d}\check{\omega}+\check{\omega}\wedge\check{\omega}$ & Field strength

$\check{F}=d\check{A}-iq\check{A}\land\check{A}$\tabularnewline
\hline 
\end{tabular*}

\caption{Constructs as applied to the various spaces associated with a point
$p\in M$ in spacetime and a vector $v$ at $p$.}
\end{table}

\begin{figure}[H]
\noindent \begin{centering}
\includegraphics[width=1\columnwidth]{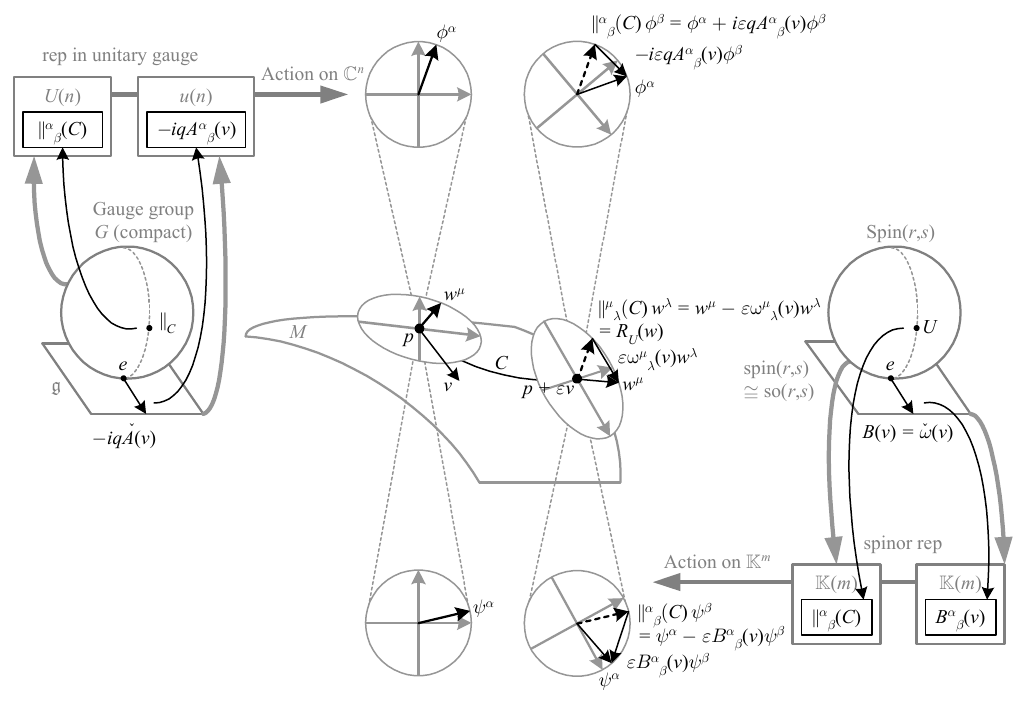}
\par\end{centering}
\caption{A matter field can be the tensor product of a complex scalar or Yang-Mills
field $\phi$ and a spinor field $\psi$. YM fields use a connection
and gauge (frame) which are independent of the spacetime manifold
frame, while spinor fields mirror the connection and changes in frame
of the spacetime manifold. YM fields are acted on by reps of the gauge
group and its Lie algebra, while spinor fields are acted on by reps
of the Spin group and its Lie algebra. In the figure we assume an
infinitesimal curve $C$ with tangent $v$, an orthonormal frame,
a spin connection, and a unitary gauge. Note that the field components
shown at $p+\varepsilon v$ are those of the field value at $p$ applied
to the frame at $p+\varepsilon v$, i.e. the top right field vector
depicted would be more precisely written $\left.\phi^{\alpha}\right|_{p}\left.e_{\alpha}\right|_{p+\varepsilon v}$;
in particular, this quantity is unrelated to the value of the field
$\left.\phi^{\alpha}\right|_{p+\varepsilon v}$. }
\end{figure}

\noindent %
\begin{framed}%
\noindent $\triangle$ Note that a Lorentz transformation on all of
flat Minkowski space, which is the setting for many treatments of
this material, induces a change of coordinate frame that is the same
Lorentz transformation on every tangent space, thus simplifying the
above picture by eliminating the need to consider parallel transport
on the curved spacetime manifold.\end{framed}%
\begin{framed}%
\noindent $\triangle$ The spinor space is an internal space, but
its changes of frame are driven by those of the spacetime manifold.
The question of whether a global change of orthonormal frame can be
mapped to globally defined elements in $\mathrm{Spin}(r,s)$ across
coordinate charts in a consistent way is resolved in Section \ref{subsec:Spinor-bundles}
in terms of fiber bundles.\end{framed}

\section{Defining bundles}

When introducing tangent spaces on a manifold $M^{n}$, the tangent
bundle is usually defined to be the set of tangent spaces at every
point within the region of a coordinate chart $U\to\mathbb{R}^{n}$,
i.e. it is defined as the cartesian product $U\times\mathbb{R}^{n}$.
Globally, one uses an atlas of charts covering $M$, with coordinate
transformations $\mathbb{R}^{n}\to\mathbb{R}^{n}$ defining how to
consider a vector field across charts. Here we want to take the same
approach to define the global version of the tangent bundle, with
analogs for frames and internal spaces.

\subsection{\label{subsec:Fiber-bundles}Fiber bundles }

In defining fiber bundles we first consider a \textbf{base space}\index{base space}
$M$ and a \textbf{bundle space}\index{bundle space} (AKA total space\index{total space},
entire space\index{entire space}) $E$, which includes a surjective
\textbf{bundle projection}\index{bundle projection} (AKA bundle submersion\index{bundle submersion},
projection map\index{projection map}) 
\begin{equation}
\pi\colon E\to M.
\end{equation}
In the special case that $M$ and $E$ are manifolds, we require the
bundle projections $\pi$ to be (infinitely) differentiable, and $E$
without any further structure is called a \textbf{fibered manifold}\index{fibered manifold}.

The space $E$ becomes a \textbf{fiber bundle}\index{fiber bundle}
(AKA fibre bundle) if each \textbf{fiber over} \textbf{\textit{x}}\index{fiber over x}
$\pi^{-1}(x)$, where $x\in M$, is homeomorphic to an \textbf{abstract
fiber}\index{standard fiber} (AKA standard fiber, typical fiber,
fiber space, fiber) $F$; specifically, we must have the analog of
an atlas, a collection of open \textbf{trivializing neighborhoods}\index{trivializing neighborhoods}
$\{U_{i}\}$ that cover $M$, each with a \textbf{local trivialization}\index{local trivialization},
a homeomorphism 
\begin{equation}
\begin{aligned}\phi_{i}\colon\pi^{-1}(U_{i}) & \to U_{i}\times F\\
p & \mapsto\left(\pi(p),f_{i}(p)\right),
\end{aligned}
\end{equation}
which in a given $\pi^{-1}(x)$ allows us to ignore the first component
and consider the last as a homeomorphism 
\begin{equation}
f_{i}\colon\pi^{-1}(x)\to F.
\end{equation}
This property of a bundle is described by calling it \textbf{locally
trivial}\index{locally trivial} (AKA a \textbf{local product space}\index{local product space}),
and if all of $M$ can be made a trivializing neighborhood, then $E$
is a \textbf{trivial bundle}\index{trivial bundle}, i.e. $E\cong M\times F$.
The topology of a non-trivial bundle can be defined via $E$ itself,
or imputed by the local trivializations. Note that if $F$ is discrete,
then $E$ is a covering space of $M$, and if $M$ is contractible,
then $E$ is trivial. If $F$ is given additional structure, $f_{i}$
must remain an isomorphism with respect to this structure.
\begin{figure}[H]
\noindent \begin{centering}
\includegraphics[width=1\columnwidth]{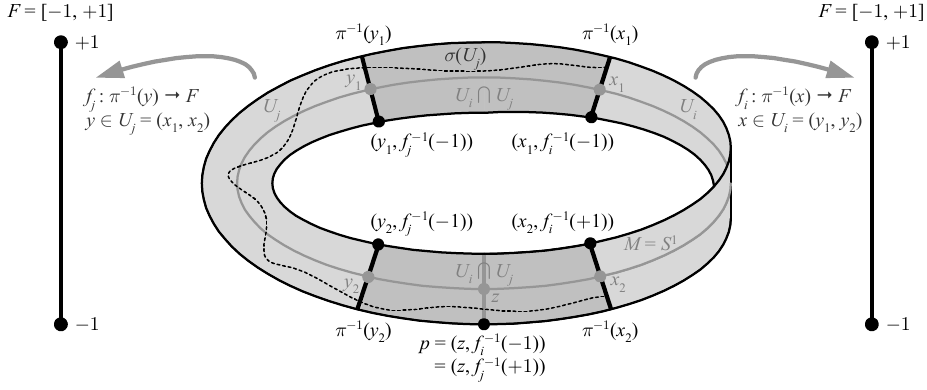}
\par\end{centering}
\caption{\label{fig:The-mobius-strip}The \textbf{Möbius strip}\index{Möbius strip}
(AKA Möbius band\index{Möbius band}) has a base space which is a
circle $M=S^{1}$, a fiber which is a line segment $F=\left[-1,+1\right]$,
and is non-trivial, since it requires at least two trivializing neighborhoods.
In the figure, the fiber over $z$ has two different descriptions
under the two local trivializations, and a local section $\sigma$
(defined below) is depicted. }
\end{figure}
\begin{framed}%
\noindent $\triangle$ Fiber bundles are denoted by various combination
of components and maps in various orders, frequently $(E,M,F)$, $(E,M,\pi)$,
or $(E,M,\pi,F)$. Other notations include $\pi\colon E\to M$ and
$F\longrightarrow E\overset{\pi}{\longrightarrow}M$ or just $F\longrightarrow E\overset{}{\longrightarrow}M$.\end{framed}%
\begin{framed}%
\noindent $\triangle$ The distinction between the fiber and the fiber
over $x$ is sometimes not made clear; it is important to remember
that the abstract fiber $F$ is not part of the bundle space $E$. \end{framed}

A \textbf{bundle map}\index{bundle map} (AKA bundle morphism\index{bundle morphism})
is a pair of maps 
\begin{equation}
\Phi_{E}\colon E\to E^{\prime}
\end{equation}
\begin{equation}
\Phi_{M}\colon M\to M^{\prime}
\end{equation}
between bundles that map fibers to fibers, i.e. 
\begin{equation}
\pi^{\prime}(\Phi_{E}(p))=\Phi_{M}(\pi(p)).
\end{equation}
Note that if $\Phi_{M}$ is the identity map, the bundles are over
the same base space $M$ and this reduces to a single map satisfying
\begin{equation}
\pi^{\prime}(\Phi(p))=\pi(p).
\end{equation}

A \textbf{section}\index{section} (AKA cross section\index{cross section})
of a fiber bundle is a continuous map 
\begin{equation}
\sigma\colon M\to E
\end{equation}
that satisfies 
\begin{equation}
\pi(\sigma(x))=x.
\end{equation}
At a point $x\in M$ a \textbf{local section}\index{local section}
always exists, being only defined in a neighborhood of $x$; however
global sections may not exist.%
\begin{framed}%
\noindent $\triangle$ It is important to remember that the base space
$M$ is not part of the bundle space $E$. In particular, since a
global section may not exist, the base space cannot in general be
viewed as being embedded in the bundle space, and even when it can
be, such an embedding is in general arbitrary. An exception is when
there is a canonical global section, for example the zero section
as depicted in the Möbius strip above (and in a vector bundle in general,
see Section \ref{subsec:Vector-bundles}).\end{framed}

\subsection{\textit{G}-bundles }

In the fiber over a point $\pi^{-1}(x)$ in the intersection of two
trivializing neighborhoods on a bundle $(E,M,F)$, we have a homeomorphism
\begin{equation}
f_{i}f_{j}^{-1}\colon F\to F.
\end{equation}
If each of these homeomorphisms is the (left) action of an element
$g_{ij}(x)\in G$, then $G$ is called the \textbf{structure group}\index{structure group}
of $E$. This action is usually required to be faithful, so that each
$g\in G$ corresponds to a distinct homeomorphism of $F$. The map
\begin{equation}
g_{ij}\colon U_{i}\cap U_{j}\to G
\end{equation}
is called a \textbf{transition function}\index{transition function};
the existence of transition functions for all overlapping charts makes
$\{U_{i}\}$ a \textbf{\textit{G}}\textbf{-atlas}\index{G-atlas}
and turns the bundle into a \textbf{\textit{G}}\textbf{-bundle}\index{G-bundle}.
Applying the action of $g_{ij}$ to an arbitrary $f_{j}(p)$ yields
\begin{equation}
f_{i}(p)=g_{ij}\left(f_{j}(p)\right).
\end{equation}
For example, the Möbius strip in the previous figure has a structure
group $G=\mathbb{Z}_{2}$, where the action of $0\in G$ is multiplication
by $+1$, and the action of $1\in G$ is multiplication by $-1$.
In the top intersection $U_{i}\cap U_{j}$, $g_{ij}=0$, so that $f_{i}$
and $f_{j}$ are identical, while in the lower intersection $g_{ij}=1$,
so that 
\begin{equation}
f_{i}(p)=g_{ij}\left(f_{j}(p)\right)=1\left(f_{j}(p)\right)=-f_{j}(p).
\end{equation}

At a point in a triple intersection $U_{i}\cap U_{j}\cap U_{k}$,
the \textbf{cocycle condition}\index{cocycle condition} 
\begin{equation}
g_{ij}g_{jk}=g_{ik}
\end{equation}
can be shown to hold, which implies 
\begin{equation}
g_{ii}=e
\end{equation}
and 
\begin{equation}
g_{ji}=g_{ij}^{-1}.
\end{equation}
Going the other direction, if we start with transition functions from
$M$ to $G$ acting on $F$ that obey the cocycle condition, then
they determine a unique $G$-bundle $E$. %
\begin{framed}%
\noindent $\triangle$ It is important to remember that the left action
of $G$ is on the abstract fiber $F$, which is not part of the entire
space $E$, and whose mappings to $E$ are dependent upon local trivializations.
A left action on $E$ itself based on these mappings cannot in general
be consistently defined, since for non-abelian $G$ it will not commute
with the transition functions.\end{framed}

A given $G$-atlas may not need all the possible homeomorphisms of
$F$ between trivializing neighborhoods, and therefore will not ``use
up'' all the possible values in $G$. If there exists trivializing
neighborhoods on a $G$-bundle whose transition functions take values
only in a subgroup $H$ of $G$, then we say the structure group $G$
is \textbf{reducible}\index{reducible structure group} to $H$. For
example, a trivial bundle's structure group is always reducible to
the trivial group consisting only of the identity element.

\subsection{\label{subsec:Principal-bundles}Principal bundles }

A \textbf{principal bundle}\index{principal bundle} (AKA principal
$G$-bundle\index{principal G-bundle}) $(P,M,\pi,G)$ has a topological
group $G$ as both abstract fiber and structure group, where $G$
acts on itself via left translation as a transition function across
trivializing neighborhoods, i.e. 
\begin{equation}
f_{i}(p)=g_{ij}f_{j}(p),
\end{equation}
where the operation of $g_{ij}$ is the group operation. Note that
the fiber over a point $\pi^{-1}(x)$ is only homeomorphic as a space
to $G$ in a given trivializing neighborhood, and so is missing a
unique identity element and is a $G$-torsor, not a group. 

A principal bundle lets us introduce a consistent right action of
$G$ on $\pi^{-1}(x)$ (as opposed to the left action on the abstract
fiber). This right action is defined by 
\begin{equation}
\begin{aligned}g(p) & \equiv f_{i}^{-1}\left(f_{i}(p)g\right)\\
\Rightarrow f_{i}\left(g(p)\right) & =f_{i}(p)g
\end{aligned}
\end{equation}
for $p\in\pi^{-1}(U_{i})$, where in an intersection of trivializing
neighborhoods $U_{i}\cap U_{j}$ we see that 
\begin{equation}
\begin{aligned}g(p) & =f_{j}^{-1}\left(f_{j}(p)g\right)\\
 & =f_{i}^{-1}f_{i}f_{j}^{-1}\left(f_{j}(p)g\right)=f_{i}^{-1}\left(g_{ij}f_{j}(p)g\right)\\
 & =f_{i}^{-1}\left(f_{i}(p)g\right)=g(p),
\end{aligned}
\end{equation}
i.e. $g(p)$ is consistently defined across trivializing neighborhoods.
Via this fiber-wise action, $G$ then has a right action on the bundle
$P$. 
\begin{figure}[H]
\noindent \begin{centering}
\includegraphics[width=1\columnwidth]{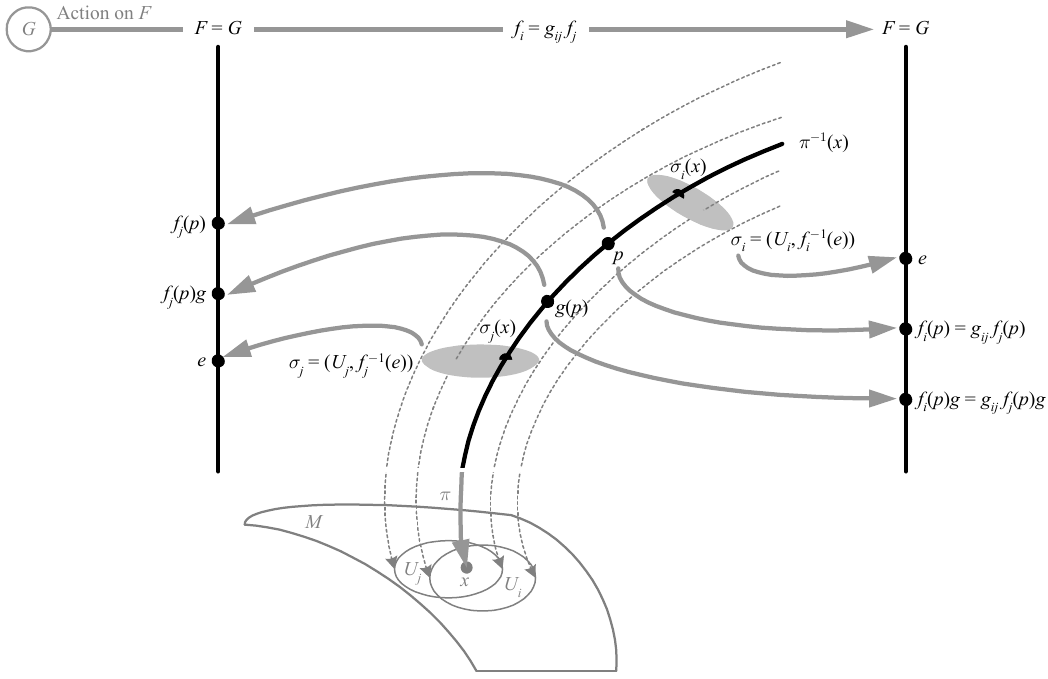}
\par\end{centering}
\caption{A principal bundle has the same group $G$ as both abstract fiber
and structure group, where $G$ acts on itself via left translation.
$G$ also has a right action on the bundle itself, which is consistent
across trivializing neighborhoods. The identity sections (defined
below) are also depicted.}
\end{figure}
\begin{framed}%
\noindent $\triangle$ It is important to remember that $M$ is not
part of $E$, and that the depiction of each fiber in the bundle $\pi^{-1}(x)\in E$
as ``hovering over'' the point $x\in M$ is only valid locally.\end{framed}%
\begin{framed}%
\noindent $\triangle$ Note that from its definition and basic group
properties, the right action of $G$ on $\pi^{-1}(x)$ is automatically
free and transitive (making $\pi^{-1}(x)$ a ``right $G$-torsor'').
An equivalent definition of a principal bundle excludes $G$ as a
structure group but includes this free and transitive right action
of $G$. Also note that the definition of the right action is equivalent
to saying that $f_{i}\colon\pi^{-1}(x)\to G$ is equivariant with
respect to the right action of $G$ on $\pi^{-1}(x)$ and the right
action of $G$ on itself.\end{framed}%
\begin{framed}%
\noindent $\triangle$ A principal bundle is sometimes defined so
that the structure group acts on itself by right translation instead
of left. In this case the action of $G$ on the bundle must be a left
action.\end{framed}%
\begin{framed}%
\noindent $\triangle$ A principal bundle can also be denoted $P(M,G)$
or $G\hookrightarrow P\overset{\pi}{\longrightarrow}M$.\end{framed}Since the right action is an intrinsic operation, a \textbf{principal
bundle map}\index{principal bundle map} between principal $G$-bundles
(e.g. a principal bundle automorphism) is required to be equivariant
with regard to it, i.e. we require 
\begin{equation}
\Phi_{E}(g(p))=g(\Phi_{E}(p)),
\end{equation}
or in juxtaposition notation, 
\begin{equation}
\Phi_{E}(pg)=\Phi_{E}(p)g.
\end{equation}
In fact, any such equivariant map is automatically a principal bundle
map, and if the base spaces are identical and unchanged by $\Phi_{E}$,
then $\Phi_{E}$ is an isomorphism. For a principal bundle map 
\begin{equation}
\Phi_{E}\colon(P^{\prime},M^{\prime},G^{\prime})\to(P,M,G)
\end{equation}
between bundles with different structure groups, we must include a
homomorphism 
\begin{equation}
\Phi_{G}\colon G^{\prime}\to G
\end{equation}
between structure groups so that the equivariance condition becomes
\begin{equation}
\Phi_{E}(g(p))=\Phi_{G}(g)(\Phi_{E}(p)),
\end{equation}
or in juxtaposition notation, 
\begin{equation}
\Phi_{E}(pg)=\Phi_{E}(p)\Phi_{G}(g).
\end{equation}
\begin{framed}%
\noindent $\triangle$ Note that the right action of a fixed $g\in G$
is thus not a principal bundle automorphism, since for non-abelian
$G$ it will not commute with another right action. \end{framed}

A principal bundle has a global section iff it is trivial. However,
within each trivializing neighborhood on a principal bundle we can
define a local \textbf{identity section}\index{identity section}
\begin{equation}
\sigma_{i}(x)\equiv f_{i}^{-1}(e),
\end{equation}
where $e$ is the identity element in $G$. In $U_{i}\cap U_{j}$,
we can then use $f_{i}(\sigma_{i})=e$ to see that the identity sections
are related by the right action of the transition function: 
\begin{equation}
\begin{aligned}g_{ij}(\sigma_{i}) & =f_{i}^{-1}\left(f_{i}(\sigma_{i})g_{ij}\right)\\
 & =f_{i}^{-1}\left(g_{ij}\right)\\
 & =f_{i}^{-1}\left(g_{ij}f_{j}(\sigma_{j})\right)\\
 & =f_{i}^{-1}\left(f_{i}(\sigma_{j})\right)\\
 & =\sigma_{j},
\end{aligned}
\end{equation}
or in juxtaposition notation,
\begin{equation}
\sigma_{j}=\sigma_{i}g_{ij}.
\end{equation}
\begin{framed}%
\noindent $\triangle$ The different actions of $G$ are a potential
source of confusion. $g_{ij}$ has a left action on the abstract fiber
of a $G$-bundle, which on a principal bundle becomes left group multiplication,
and also has a right action on the bundle itself that relates the
elements in the identity section. \end{framed}

If $G$ is a closed subgroup of a Lie group $P$ (and thus also a
Lie group by Cartan's theorem), then $(P,P/G,G)$ is a principal $G$-bundle
with base space the (left) coset space $P/G$. The right action of
$G$ on the entire space $P$ is just right translation. 

\section{Generalizing tangent spaces}

In this section we use matrix notation to reduce clutter, remembering
that bases are row vectors and are acted on by matrices from the right.
We retain index notation when acting on vector components to avoid
confusion with operations on intrinsic vectors.

\subsection{\label{subsec:Associated-bundles}Associated bundles }

If two $G$-bundles $(E,M,F)$ and $(E^{\prime},M,F^{\prime})$, with
the same base space and structure group, also share the same trivializing
neighborhoods and transition functions, then they are each called
an \textbf{associated bundle}\index{associated bundle} with regard
to the other. It is possible to construct (up to isomorphism) a unique
principal $G$-bundle associated to a given $G$-bundle; going in
the other direction, given a principal $G$-bundle and a left action
of $G$ on a fiber $F$, we can construct a unique associated $G$-bundle
with fiber $F$. In particular, given a principal bundle $(P,M,G)$,
the rep of $G$ on itself by inner automorphisms defines an associated
bundle $(\mathrm{Inn}P,M,G)$, and the adjoint rep of $G$ on $\mathfrak{g}$
defines an associated bundle $(\mathrm{Ad}P,M,\mathfrak{g})$. If
$G$ has a linear rep on a vector space $\mathbb{K}^{n}$, this rep
defines an associated bundle $(E,M,\mathbb{K}^{n})$, which we explore
next. 
\begin{figure}[H]
\noindent \begin{centering}
\includegraphics[width=1\columnwidth]{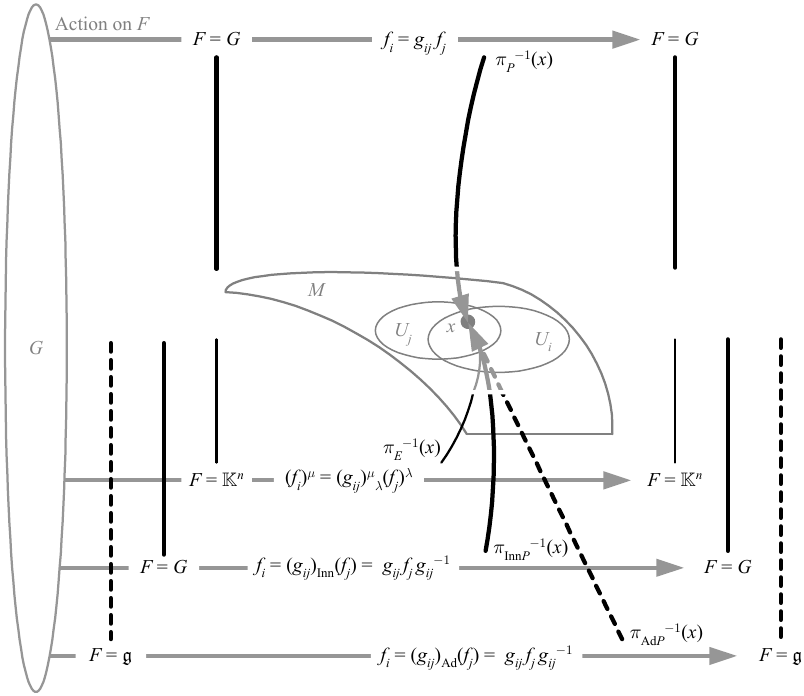}
\par\end{centering}
\caption{Given a principal bundle, we can construct an associated bundle for
the action of $G$ on a vector space $\mathbb{K}^{n}$ by a linear
rep, on itself by inner automorphisms, and on its Lie algebra $\mathfrak{g}$
by the adjoint rep. The action of the structure group is shown in
general and for the case in which $G$ is a matrix group, with matrix
multiplication denoted as juxtaposition. Although denoted identically,
the $f_{i}$ are those corresponding to each bundle. }
\end{figure}
\begin{framed}%
\noindent $\triangle$ The $G$-bundle $E$ with fiber $F$ associated
to a principal bundle $P$ is sometimes written 
\begin{equation}
E=P\times_{G}F\equiv(P\times F)/G,
\end{equation}
where the quotient space collapses all points in the product space
which are related by the right action of some $g\in G$ on $P$ and
the right action of $g^{-1}$ on $F$. \end{framed}

\subsection{\label{subsec:Vector-bundles}Vector bundles }

A \textbf{vector bundle}\index{vector bundle} $(E,M,\pi,\mathbb{K}^{n})$
has a vector space fiber $\mathbb{K}^{n}$ (assumed here to be $\mathbb{R}^{n}$
or $\mathbb{C}^{n}$) and a structure group that is linear ($G\subseteq GL(n,\mathbb{K})$)
and therefore acts as a matrix across trivializing neighborhoods,
i.e. 
\begin{equation}
f_{i}(p)=g_{ij}f_{j}(p),
\end{equation}
where the operation of $g_{ij}$ is now matrix multiplication on the
vector components $f_{j}(p)\in\mathbb{K}^{n}$. If we view $V_{x}\equiv\pi^{-1}(x)$
as an internal space on $M$ with intrinsic vector elements $v$,
the linear map $f_{i}\colon\pi^{-1}(x)\to\mathbb{K}^{n}$ is equivalent
to choosing a basis $e_{i\mu}$ to get vector components, i.e. 
\begin{equation}
f_{i}(v)=v_{i}^{\mu},
\end{equation}
where 
\begin{equation}
v_{i}^{\mu}e_{i\mu}=v
\end{equation}
and latin letters are labels while greek letters are the usual indices
for vectors and labels for bases. The action of the structure group
can then be written 
\begin{equation}
v_{i}^{\mu}=(g_{ij})^{\mu}{}_{\lambda}v_{j}^{\lambda},
\end{equation}
which is equivalent to a change of basis 
\begin{equation}
e_{i\mu}=(g_{ij}^{-1})^{\lambda}{}_{\mu}e_{j\lambda},
\end{equation}
or as matrix multiplication on basis row vectors 
\begin{equation}
e_{j}=e_{i}g_{ij},
\end{equation}
so that the action of $g_{ij}(x)$ in $U_{i}\cap U_{j}$ is equivalent
to a change of frame or gauge transformation from $e_{i}$ to $e_{j}$,
which is equivalent to a transformation of internal space vector components
in the opposite direction.%
\begin{framed}%
\noindent $\triangle$ The frame is not a part of the vector bundle,
it is a way of viewing the local trivializations; therefore the view
of $g_{ij}(x)$ as effecting a change of basis should not be confused
with a group action on either $\pi^{-1}(x)$ or $E$. As the structure
group of $E$, the action of $G$ is on the fiber $\mathbb{K}^{n}$,
which is not part of $E$.\end{framed}

If the structure group of a vector bundle is reducible to $GL(n,\mathbb{K})^{e}$,
then it is called an \textbf{orientable bundle}\index{orientable bundle};
all complex vector bundles are orientable, so orientability usually
refers to real vector bundles. The tangent bundle of $M$ (formally
defined in Section \ref{subsec:The-tangent-bundle-and-solder-form})
is then orientable iff $M$ is orientable. On a pseudo-Riemannian
manifold $M$, the structure group is reducible to $O(r,s)$, and
if $M$ is orientable then it is reducible to $SO(r,s)$; if the structure
group can be further reduced to $SO(r,s)^{e}$, then $M$ and its
tangent bundle are called \textbf{time and space orientable}\index{time and space orientable}\index{orientable!time and space}.
Note that this additional distinction is dependent only upon the metric,
and two metrics on the same manifold can have different time and space
orientabilities. %
\begin{framed}%
\noindent $\triangle$ The orientability of a vector bundle as a bundle
is different than its orientability as a manifold itself; therefore
it is important to understand which version of orientability is being
referred to. In particular, the tangent bundle of $M$ is always orientable
as a manifold, but it is orientable as a bundle only if $M$ is.\end{framed}

A gauge transformation on a vector bundle is a smoothly defined linear
transformation of the basis inferred by the components due to local
trivializations at each point, i.e. 
\begin{equation}
e_{i\mu}^{\prime}=(\gamma_{i}^{-1})^{\lambda}{}_{\mu}e_{i\lambda},
\end{equation}
which is equivalent to new local trivializations where
\begin{equation}
\left(v_{i}^{\mu}\right)^{\prime}=(\gamma_{i})^{\mu}{}_{\lambda}v_{i}^{\lambda},
\end{equation}
giving us new transition functions
\begin{equation}
g_{ij}^{\prime}=\gamma_{i}g_{ij}\gamma_{j}^{-1},
\end{equation}
where we have suppressed indices for pure matrix relationships. Thus
the gauge group is the same as the structure group, and a gauge transformation
$\gamma_{i}^{-1}$ is equivalent to the transition function $g_{i^{\prime}i}$
from $U_{i}$ to $U_{i}^{\prime}$, the same neighborhood with a different
local trivialization. 
\begin{figure}[H]
\noindent \begin{centering}
\includegraphics[width=1\columnwidth]{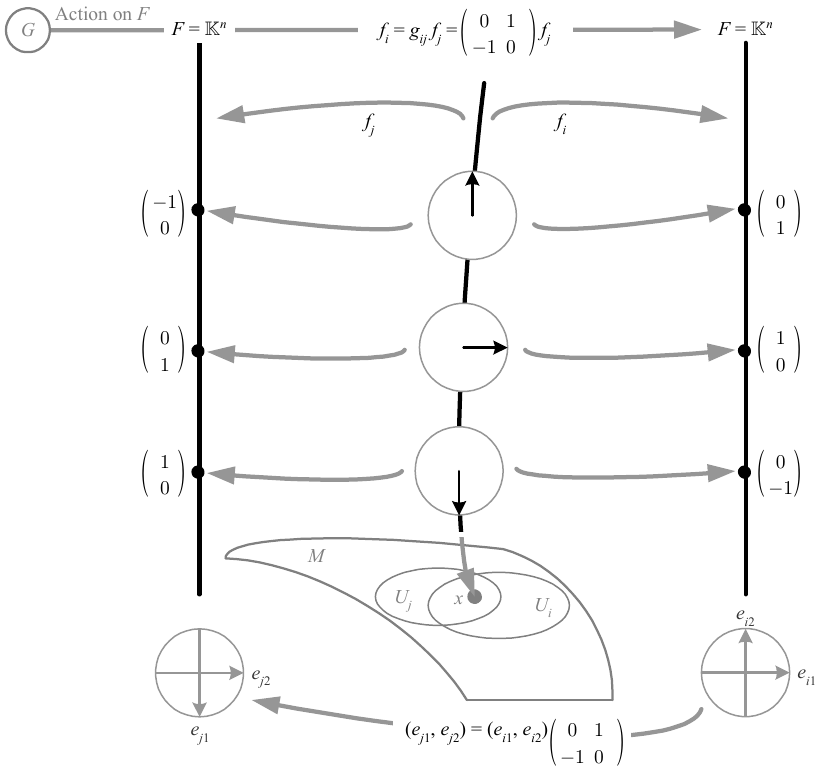}
\par\end{centering}
\caption{The elements of the fiber over $x$ in a vector bundle can be viewed
as abstract vectors in an internal space, with the local trivialization
acting as a choice of basis from which the components of these vectors
can be calculated. The structure group then acts as a matrix transformation
between vector components, and between bases in the opposite direction.
A gauge transformation is also a new choice of basis, and so can be
handled similarly.}
\end{figure}

A vector bundle always has global sections (e.g. the zero vector in
the fiber over each point). A vector bundle with fiber $\mathbb{R}$
is called a \textbf{line bundle}\index{line bundle}. 

\subsection{Frame bundles }

Given a vector bundle $(E,M,\mathbb{K}^{n})$, the \textbf{frame bundle}\index{frame bundle}
of $E$ is the principal $GL(n,\mathbb{K})$-bundle associated to
$E$, and is denoted 
\begin{equation}
F(E)\equiv(F(E),M,\pi,GL(n,\mathbb{K})).
\end{equation}
The elements $p\in\pi^{-1}(x)$ are viewed as ordered bases of the
internal space $V_{x}\cong\mathbb{K}^{n}$, which we denote 
\begin{equation}
p\equiv e_{p},
\end{equation}
or $e_{p\mu}$ if operated on by a matrix in index notation. Each
trivializing neighborhood $U_{i}$ is associated with a fixed frame
$e_{i}$, which we take from the local trivializations in the vector
bundle $E$, letting us define 
\begin{equation}
f_{i}\colon\pi^{-1}(x)\to GL(n,\mathbb{K})
\end{equation}
by the matrix relation 
\begin{equation}
e_{p}=e_{i}f_{i}(p).
\end{equation}
In other words $f_{i}(p)$ is the matrix that transforms (as row vectors)
the fixed basis $e_{i}$ into the basis element $e_{p}$ of $F(E)$;
in particular, the identity section is 
\begin{equation}
\sigma_{i}=f_{i}^{-1}(I)=e_{i},
\end{equation}
where $I$ is the identity matrix. If we again write vector components
in these bases as $v_{i}^{\mu}e_{i\mu}=v_{p}^{\mu}e_{p\mu}=v$, then
we have 
\begin{equation}
v_{i}^{\mu}=f_{i}(p)^{\mu}{}_{\lambda}v_{p}^{\lambda}.
\end{equation}

The left action of $g_{ij}$ is defined by $f_{i}(p)=g_{ij}f_{j}(p)$,
and applying both sides to vector components $v_{p}^{\mu}$ we get
\begin{equation}
v_{i}^{\mu}=(g_{ij})^{\mu}{}_{\lambda}v_{j}^{\lambda},
\end{equation}
the same transition functions as in $E$. The transition functions
can be viewed as changes of frame $e_{j}=e_{i}g_{ij}$, or gauge transformations,
between the identity sections of $F(E)$ in $U_{i}\cap U_{j}$, i.e.
this can be written as a matrix relation 
\begin{equation}
\sigma_{j}=\sigma_{i}g_{ij},
\end{equation}
which as we see next is the usual right action of the transition functions
on identity sections.
\begin{figure}[H]
\noindent \begin{centering}
\includegraphics[width=1\columnwidth]{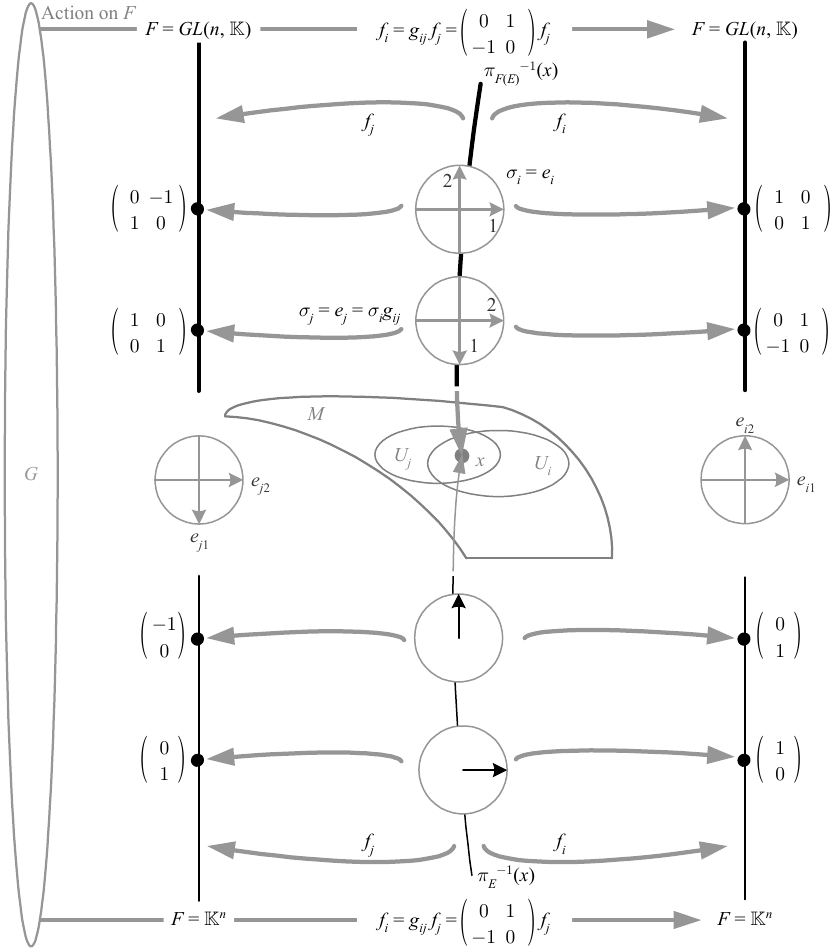}
\par\end{centering}
\caption{Given a vector bundle $E$, we can construct an associated frame bundle
$F(E)$. The elements of the fiber over $x$ in the frame bundle can
be viewed as bases for the internal space, with the local trivialization
acting as a choice of a fixed basis against which linear transformations
generate these bases. These fixed bases are the same as those chosen
in the corresponding local trivialization on the vector bundle, and
are acted on by the same transition functions. Although denoted identically,
the $f_{i}$ are those corresponding to each bundle. }
\end{figure}
\begin{framed}%
\noindent $\triangle$ Unlike with $E$, the frame is in fact part
of the bundle $F(E)$, but vectors and vector components are not.
The left action of $g_{ij}$ on the abstract fiber $GL(n,\mathbb{K})$
is equivalent to a transformation in the opposite direction from the
fixed frame in $U_{i}$ to the fixed frame in $U_{j}$, which is a
right action on the identity sections from $\sigma_{i}=e_{i}$ to
$\sigma_{j}=e_{j}$.\end{framed}%
\begin{framed}%
\noindent $\triangle$ It is important to remember that the elements
of $\pi^{-1}(x)$ in $F(E)$ are bases of the vector space $V_{x}$,
and in a given trivializing neighborhood it is only the basis in the
identity section that is identified with the basis underlying the
vector components in the same trivializing neighborhood of $E$. \end{framed}

The right action of $g\in GL(n,\mathbb{K})$ on $\pi^{-1}(x)$ is
defined by $f_{i}(g(p))=f_{i}(p)g$, and applying both sides to $e_{i}$
from the right and using $e_{p}=e_{i}f_{i}(p)$ we immediately obtain
\begin{equation}
e_{g(p)}=e_{p}g,
\end{equation}
so that the right action of the matrix $g$ is literally matrix multiplication
from the right on the basis row vector $p=e_{p\mu}$. Alternative
ways of writing this relation include 
\begin{equation}
\begin{aligned}e_{g(p)\mu} & =e_{p\mu}g^{\mu}{}_{\lambda},\\
g(p) & =pg.
\end{aligned}
\end{equation}
In particular, if $f_{i}(p)=g$ then we have 
\begin{equation}
p=e_{p}=e_{i}g=g(e_{i})=e_{g(e_{i})}.
\end{equation}
\begin{framed}%
\noindent \sun{} Note that since the right action on $\pi^{-1}(x)$
is by a fixed matrix, it acts as a transformation relative to each
$e_{p}$, not as a transformation on the internal space $V_{x}$ in
which all of the bases in $\pi^{-1}(x)$ live. As a concrete example,
if $g^{0}{}_{0}=1$ and $g^{\lambda\neq0}{}_{0}=0$, then $e_{g(p)0}=e_{p0}$,
meaning that the transformation $p\mapsto g(p)$ leaves first vector
of all bases in $\pi^{-1}(x)$ unaffected, regardless of that vector's
direction. This behavior contrasts with that of a transformation on
$V_{x}$ itself, which as we will see in the next section is a gauge
transformation.\end{framed}

\subsection{\label{subsec:Gauge-transformations-on-frame-bundles}Gauge transformations
on frame bundles }

Recall that a gauge transformation on a vector bundle $E$ is an active
transformation of the bases underlying the components defining a local
trivialization, which is equivalent to a new set of local trivializations
and transition functions (and is not a transformation on the space
$E$ itself). On the frame bundle $F(E)$, we perform the same basis
change for the fixed frames associated with each trivializing neighborhood
\begin{equation}
e_{i}^{\prime}=e_{i}\gamma_{i}^{-1},
\end{equation}
which also defines the new identity sections, and is equivalent to
new local trivializations where
\begin{equation}
f_{i}^{\prime}(p)=\gamma_{i}f_{i}(p),
\end{equation}
giving us new transition functions
\begin{equation}
g_{ij}^{\prime}=\gamma_{i}g_{ij}\gamma_{j}^{-1},
\end{equation}
which are the same as those in the associated vector bundle $E$.
We will call this transformation a \textbf{neighborhood-wise gauge
transformation}\index{neighborhood-wise gauge transformation}\index{gauge transformation!neighborhood-wise}. 

An alternative (and more common) way to view gauge transformations
on $F(E)$ is to transform the actual bases in $\pi^{-1}(x)$ via
a bundle automorphism 
\begin{equation}
p^{\prime}\equiv\gamma^{-1}(p),
\end{equation}
and then change the fixed bases in each trivializing neighborhood
to 
\begin{equation}
\begin{aligned}e_{i}^{\prime} & =\gamma^{-1}(e_{i})\\
 & \equiv e_{i}\gamma_{i}^{-1}
\end{aligned}
\end{equation}
in order to leave the maps $f_{i}(p)$ the same (which also leaves
the identity sections and transition functions the same). This immediately
implies a constraint on the basis changes in $U_{i}\cap U_{j}$: since
$g_{ij}^{\prime}=\gamma_{i}g_{ij}\gamma_{j}^{-1}$, requiring constant
$g_{ij}$ means we must have
\begin{equation}
\gamma_{i}^{-1}=g_{ij}\gamma_{j}^{-1}g_{ij}^{-1}.
\end{equation}
We will call this transformation an \textbf{automorphism gauge transformation}\index{automorphism gauge transformation}\index{gauge transformation!automorphism}.%
\begin{framed}%
\noindent $\triangle$ Note that this constraint means that automorphism
gauge transformations are a subset of neighborhood-wise gauge transformations,
which allow arbitrary changes of frame in every trivializing neighborhood.
Also note that for automorphism gauge transformations, the matrices
$\gamma_{i}^{-1}$ (and therefore the new identity section elements
$e_{i}^{\prime}$) are determined by the automorphism $\gamma^{-1}$,
while neighborhood-wise gauge transformations are defined by arbitrary
matrices $\gamma_{i}^{-1}$ in each neighborhood which are not necessarily
consistent in $U_{i}\cap U_{j}$.\end{framed}%
\begin{framed}%
\noindent \sun{} As with the associated vector bundle, for either
type of gauge transformation the gauge group is the same as the structure
group, and a gauge transformation $\gamma_{i}^{-1}$ is equivalent
to the transition function $g_{i^{\prime}i}$ from $U_{i}$ to $U_{i}^{\prime}$,
the same neighborhood with a different local trivialization.\end{framed}We now define the matrices $\gamma_{p}^{-1}$ to be those which result
from the transformation $\gamma^{-1}(p)$ on the rest of $\pi^{-1}(x)$,
i.e. 
\begin{equation}
e_{p}^{\prime}\equiv e_{p}\gamma_{p}^{-1}.
\end{equation}
Note that $\gamma_{p}^{-1}$ is determined by $\gamma_{i}^{-1}$:
since we require that $f_{i}^{\prime}=f_{i}$, we have 
\begin{equation}
\begin{aligned}e_{i}^{\prime}f_{i}(p) & =e_{p}^{\prime}\\
\Rightarrow e_{i}\gamma_{i}^{-1}f_{i}(p) & =e_{p}\gamma_{p}^{-1}\\
 & =e_{i}f_{i}(p)\gamma_{p}^{-1}\\
\Rightarrow\gamma_{p}^{-1} & =f_{i}(p)^{-1}\gamma_{i}^{-1}f_{i}(p),
\end{aligned}
\end{equation}
or more generally, using the definition of a right action $f_{i}(g(p))=f_{i}(p)g$
we get
\begin{equation}
\gamma_{g(p)}^{-1}=g^{-1}\gamma_{p}^{-1}g.
\end{equation}

\noindent %
\begin{framed}%
\noindent $\triangle$ It is important to remember that the matrices
$\gamma_{i}^{-1}$ are dependent upon the local trivialization (since
they are defined as the matrix acting on the element $e_{i}\in\pi^{-1}(x)$
for $x\in U_{i}$), but the matrices $\gamma_{p}^{-1}$ are independent
of the local trivialization, and are the action of the automorphism
$\gamma^{-1}$ on the basis $e_{p}$. \end{framed}
\begin{figure}[H]
\noindent \begin{centering}
\includegraphics[width=1\columnwidth]{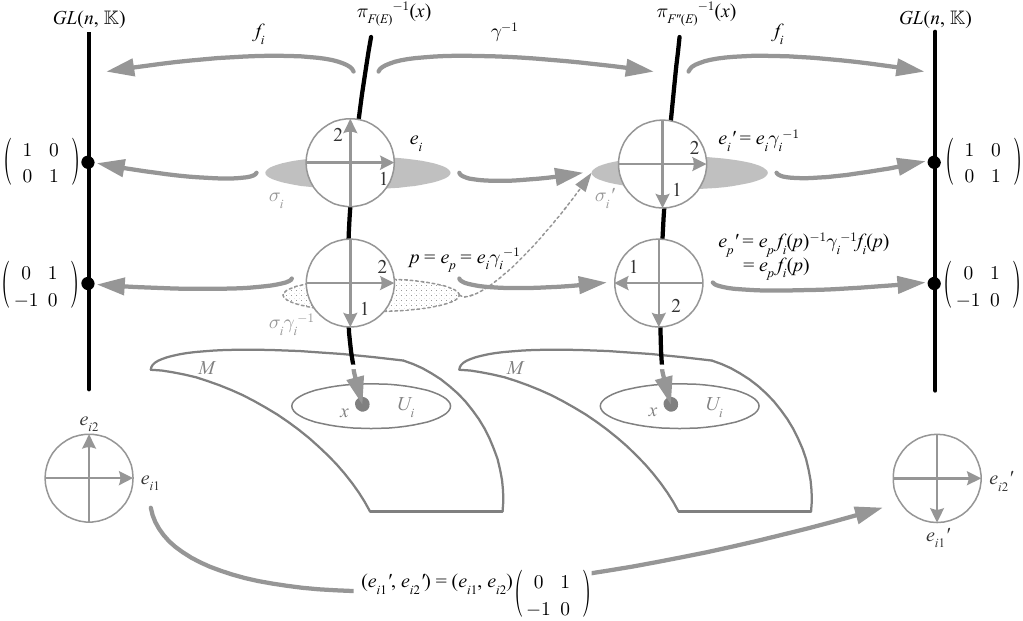}
\par\end{centering}
\caption{An automorphism gauge transformation on $F(E)$ transforms the actual
elements of the fiber over $x$, including the identity section elements
corresponding to the fixed bases in each local trivialization, thus
leaving the local trivializations unchanged.}
\end{figure}
\begin{framed}%
\noindent \sun{} This result can be understood as $\gamma^{-1}$ being
a transformation on the internal space $V_{x}$ itself, applied to
all the elements of $\pi^{-1}(x)$, each of which is a basis of $V_{x}$.
For example, in the figure above, $\gamma^{-1}$ rotates all bases
clockwise by $\pi/2$. To see why this is so, note that the matrix
in the transformation $\left(v_{i}^{\mu}\right)^{\prime}=(\gamma_{i})^{\mu}{}_{\lambda}v_{i}^{\lambda}$
has components which are those of $\gamma_{i}\in GL(V_{x})$ in the
basis $e_{i\mu}$. Therefore in a different basis $e_{p\mu}\in\pi^{-1}(x)$
we must apply a different matrix $\left(v_{p}^{\mu}\right)^{\prime}=(\gamma_{p})^{\mu}{}_{\lambda}v_{p}^{\lambda}$
which reflects the change of basis $e_{p\mu}=f_{i}(p)^{\lambda}{}_{\mu}e_{i\lambda}$
via a similarity transformation 
\begin{equation}
\begin{aligned}\gamma_{p} & =f_{i}(p)^{-1}\gamma_{i}f_{i}(p)\\
\Rightarrow\gamma_{p}^{-1} & =f_{i}(p)^{-1}\gamma_{i}^{-1}f_{i}(p).
\end{aligned}
\end{equation}

Viewed as a transformation on $V_{x}$, $\gamma^{-1}$ will then commute
with any fixed matrix applied to the bases, which as we saw is the
right action; as we see next, this corresponds to the equivariance
of $\gamma^{-1}$ required by it being a bundle automorphism.\end{framed}We now check that $\gamma^{-1}$ is a bundle automorphism with respect
to the right action of $G$, i.e. that $\gamma^{-1}\left(g(p)\right)=g\left(\gamma^{-1}(p)\right)$:

\begin{equation}
\begin{aligned}\gamma^{-1}\left(g(p)\right) & =e_{g(p)}\gamma_{g(p)}^{-1}\\
 & =e_{g(p)}g^{-1}\gamma_{p}^{-1}g\\
 & =e_{p}\gamma_{p}^{-1}g\\
 & =e_{\gamma^{-1}(p)}g\\
 & =g\left(\gamma^{-1}(p)\right)
\end{aligned}
\end{equation}
\begin{framed}%
\noindent $\triangle$ A potential source of confusion is that a local
gauge transformation (different at different points) can be defined
globally on $F(E)$; meanwhile, a global gauge transformation (the
same matrix $\gamma_{i}^{-1}$ at every point) can only be defined
locally (unless $F(E)$ is trivial).\end{framed}

Consider the associated bundle to $F(E)$ with fiber $GL(\mathbb{K}^{n})$,
where the local trivialization of the fiber over $x$ is defined to
be the possible automorphism gauge transformations $\gamma_{i}^{-1}$
on the identity section element over $x$ in the trivializing neighborhood
$U_{i}$. Then recalling that $\gamma_{i}^{-1}=g_{ij}\gamma_{j}^{-1}g_{ij}^{-1}$,
we see that the action of the structure group on the fiber is by inner
automorphism. Since the values of $\gamma^{-1}$ on $F(E)$ are determined
by those in the identity section, we can thus view automorphism gauge
transformations as sections of the associated bundle $(\mathrm{Inn}F(E),M,GL(\mathbb{K}^{n}))$.

\subsection{Smooth bundles and jets}

Nothing we have done so far has required the spaces of a fiber bundle
to be manifolds; if they are, then we require the bundle projections
$\pi$ to be (infinitely) differentiable and $\pi^{-1}(x)$ to be
diffeomorphic to $F$, resulting in a \textbf{smooth bundle}\index{smooth bundle}.
A \textbf{smooth }\textbf{\textit{G}}\textbf{-bundle}\index{smooth G-bundle}
then has a structure group $G$ which is a Lie group, and whose elements
correspond to diffeomorphisms of $F$. 

If we consider a local section $\sigma$ of a smooth fiber bundle
$(E,M,\pi,F)$ with $\sigma(x)=p$, the equivalence class of all local
sections that have both $\sigma(x)=p$ and also the same tangent space
$T_{p}\sigma$ is called the \textbf{jet}\index{jet} $j_{p}\sigma$
with \textbf{representative}\index{jet!representative} $\sigma$.
We can also require that further derivatives of the section match
the representative, in which case the order of matching derivatives
defines the \textbf{order}\index{jet!order} of the jet, which is
also called a \textbf{\textit{k}}\textbf{-jet}\index{k-jet} so that
the above definition would be that of a 1-jet. $x$ is called the
\textbf{source}\index{jet!source} of the jet and $p$ is called its
\textbf{target}\index{jet!target}. With some work to transition between
local sections, one can then form a \textbf{jet manifold}\index{jet manifold}
by considering jets with all sources and representative sections,
which becomes a \textbf{jet bundle}\index{jet bundle} by considering
jets to be fibers over their source. 
\begin{figure}[H]
\noindent \begin{centering}
\includegraphics[width=0.6\columnwidth]{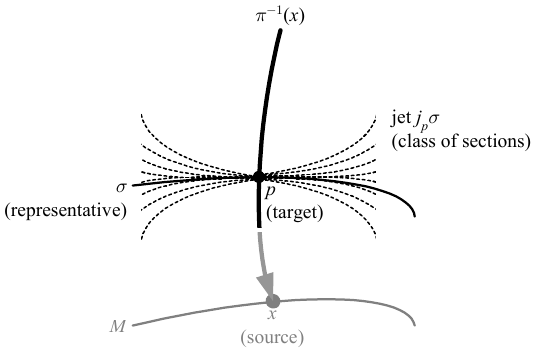}
\par\end{centering}
\caption{A jet with representative $\sigma$, source $x$, and target $p$. }
\end{figure}

\subsection{\label{subsec:Vertical-tangents-and-horizontal-equivariant-forms}Vertical
tangents and horizontal equivariant forms}

A smooth bundle $(E,M,\pi)$ is a manifold itself, and thus has tangent
vectors. A tangent vector $v$ at $p\in E$ is called a \textbf{vertical
tangent}\index{vertical tangent} if 
\begin{equation}
\mathrm{d}\pi(v)=0,
\end{equation}
i.e. if it is tangent to the fiber over $x$ where $\pi(p)=x$, so
the projection down to the base space vanishes. The \textbf{vertical
tangent space}\index{vertical tangent space} $V_{p}$ is then the
subspace of the tangent space $T_{p}$ at $p$ consisting of vertical
tangents, and viewing the vertical tangent spaces as fibers over $E$
we can form the \textbf{vertical bundle}\index{vertical bundle} $(VE,E,\pi_{V})$,
which is a subbundle of $TE$. We can also consider differential forms
on a smooth bundle, which take arguments that are tangent vectors
on $E$. A form is called a \textbf{horizontal form }\index{horizontal form}
if it vanishes whenever any of its arguments are vertical. 

On a smooth principal bundle $(P,M,G)$, we have a consistent right
action 
\begin{equation}
\rho\colon G\to\mathrm{Diff}(P),
\end{equation}
and the corresponding Lie algebra action 
\begin{equation}
\mathrm{d}\rho\colon\mathfrak{g}\to\mathrm{vect}(P)
\end{equation}
is then a Lie algebra homomorphism. The fundamental vector fields
corresponding to elements of $\mathfrak{g}$ are vertical tangent
fields; in fact, at a point $p$, $\mathrm{d}\rho\left|_{p}\right.$
is a vector space isomorphism from $\mathfrak{g}$ to $V_{p}$: 
\begin{equation}
\mathrm{d}\rho\left|_{p}\right.\colon\mathfrak{g}\overset{\cong}{\to}V_{p}
\end{equation}
In addition, the right action 
\begin{equation}
g\colon P\to P
\end{equation}
of a given element $g$ corresponds to a right action 
\begin{equation}
\mathrm{d}g\colon TP\to TP,
\end{equation}
which maps tangent vectors on $P$ via 
\[
\mathrm{d}g(v)\colon T_{p}P\to T_{g(p)}P.
\]
This map is an automorphism of $TP$ restricted to $\pi_{P}^{-1}(x)$,
which we denote $T_{\pi^{-1}(x)}P$, and it preserves vertical tangent
vectors. We can then consider the pullback 
\begin{equation}
g^{*}\varphi(v_{1},\ldots,v_{k})=\varphi(\mathrm{d}g(v_{1}),\ldots,\mathrm{d}g(v_{k}))
\end{equation}
as a right action on the space $\Lambda^{k}P$ of $k$-forms on $P$.

If we have a bundle $(E,M,\pi_{E},F)$ associated to $(P,M,\pi_{P},G)$,
we can define an $F$-valued form $\varphi_{P}$, which can be viewed
on each $\pi_{P}^{-1}(x)$ as a mapping 
\begin{equation}
\varphi_{P}\colon T_{\pi^{-1}(x)}P\otimes\cdots\otimes T_{\pi^{-1}(x)}P\to F\times\pi_{P}^{-1}(x),
\end{equation}
where $g\in G$ has a right action $\mathrm{d}g$ on $T_{\pi^{-1}(x)}P$
and a left action $g$ on the abstract fiber $F$ of $E$. The form
$\varphi_{P}$ is called an \textbf{equivariant form}\index{equivariant form}
if this mapping is equivariant with respect to these actions, i.e.
if 
\begin{equation}
g^{*}\varphi_{P}=g^{-1}\left(\varphi_{P}\right).
\end{equation}
If $\varphi_{P}$ is also horizontal, then it is called a \textbf{horizontal
equivariant form}\index{horizontal equivariant form} (AKA basic form\index{basic form},
tensorial form\index{tensorial form}). If we pull back a horizontal
equivariant form to the base space $M$ using the identity sections,
we get forms 
\begin{equation}
\varphi_{i}\equiv\sigma_{i}^{*}\varphi_{P}
\end{equation}
on each $U_{i}\subset M$. Using the identity section relation $\sigma_{i}=g_{ij}^{-1}(\sigma_{j})$
and the pullback composition property $\left(g(h)\right)^{*}\varphi=h^{*}\left(g^{*}\varphi\right)$,
we see that the values of these forms satisfy 
\begin{equation}
\begin{aligned}\varphi_{i} & =\left(g_{ij}^{-1}(\sigma_{j})\right)^{*}\varphi_{P}\\
 & =\sigma_{j}^{*}\left(\left(g_{ij}^{-1}\right)^{*}\varphi_{P}\right)\\
 & =\sigma_{j}^{*}\left(g_{ij}\left(\varphi_{P}\right)\right)\\
 & =g_{ij}\left(\varphi_{j}\right),
\end{aligned}
\end{equation}
where in the third line $g_{ij}$ is acting on the value of $\varphi_{P}$.
This means that at a point $x$ in $U_{i}\cap U_{j}$, the values
of $\varphi_{i}$ and $\varphi_{j}$ in the abstract fiber $F$ correspond
to a single point in $\pi_{E}^{-1}(x)\in E$, so that the union $\bigcup\varphi_{i}$
can be viewed as comprising a single $E$-valued form $\varphi$ on
$M$. Such a form is sometimes called a \textbf{section-valued form}\index{section-valued form},
since for fixed argument vector fields its value on $M$ is a section
of $E$. It can be shown that the correspondence between the $E$-valued
forms $\varphi$ on $M$ and the horizontal equivariant $F$-valued
forms on $P$ is one-to-one. Equivariant $F$-valued 0-forms on $P$
are automatically horizontal (since one cannot pass in a vertical
argument), and are thus one-to-one with sections on $E$. 
\begin{figure}[H]
\noindent \begin{centering}
\includegraphics[width=0.85\columnwidth]{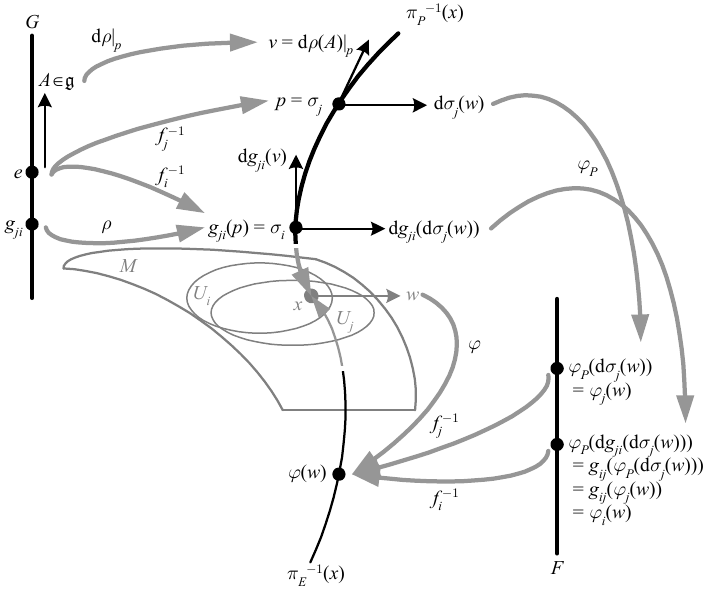}
\par\end{centering}
\caption{The differential of the right action of $G$ on $\pi_{P}^{-1}(x)\in P$
creates an isomorphism to the vertical tangent space $\mathfrak{g}\protect\cong V_{p}$.
A horizontal equivariant form $\varphi_{P}$ on $P$ maps non-vertical
vectors to the abstract fiber $F$ of an associated bundle, and pulling
back by the identity sections yields an $E$-valued form $\varphi$
on $M$. Although denoted identically, the $f_{i}$ are those corresponding
to each bundle. }
\end{figure}

On the frame bundle $(P,M,\pi_{P},GL(n,\mathbb{K}))$ associated with
a vector bundle $(E,M,\pi_{E},\mathbb{K}^{n})$, a $\mathbb{K}^{n}$-valued
form $\vec{\varphi}_{P}$ is then equivariant if 
\begin{equation}
g^{*}\vec{\varphi}_{P}=\check{g}^{-1}\vec{\varphi}_{P},
\end{equation}
where $\check{g}^{-1}$ is a matrix-valued 0-form on $P$ operating
on the $\mathbb{K}^{n}$-valued form $\vec{\varphi}_{P}$. The pullback
of a horizontal equivariant form on $P$ to the base space $M$ using
the identity sections satisfies 
\begin{equation}
\vec{\varphi}_{i}=\check{g}_{ij}\vec{\varphi}_{j},
\end{equation}
where $\check{g}_{ij}$ is now a matrix-valued 0-form on $M$. At
a point $x$ in $U_{i}\cap U_{j}$, the values of $\vec{\varphi}_{i}$
and $\vec{\varphi}_{j}$ in the abstract fiber $\mathbb{K}^{n}$ correspond
to a single abstract vector in $V_{x}=\pi_{E}^{-1}(x)\in E$, so that
the union $\bigcup\vec{\varphi}_{i}$ can be viewed as comprising
a single $V$-valued form $\vec{\varphi}$ on $M$. Thus an equivariant
$\mathbb{K}^{n}$-valued 0-form on $P$ is a matter field on $M$.
\begin{framed}%
\noindent \sun{} This correspondence can be viewed as follows. The
right action of $g$ on $P$ is a transformation on bases, so that
the equivalent transformation of vector components is $g^{-1}$. The
left action of $g^{-1}$ on the abstract fiber of $E$ is also a transformation
of vector components. Thus the equivariant property can be viewed
as ``keeping the same value when changing basis on both bundles,''
so that the values of $\vec{\varphi}_{P}$ on $\pi_{P}^{-1}(x)\in P$
correspond to a single point in $\pi_{E}^{-1}(x)\in E$, i.e a single
abstract vector over $M$. In other words, $\vec{\varphi}\in T_{x}M$
is determined by the value of $\vec{\varphi}_{P}$ at a single point
in $\pi_{P}^{-1}(x)\in P$. The horizontal requirement means we do
not consider forms which take non-zero values given argument vectors
which project down to a zero vector on $M$. \end{framed}

Under an automorphism gauge transformation, the transformation of
a horizontal equivariant form on the frame bundle $P$ is defined
by the pullback of the automorphism 
\begin{equation}
\vec{\varphi}_{P}^{\prime}\equiv\left(\gamma^{-1}\right)^{*}\vec{\varphi}_{P}.
\end{equation}
The automorphism does not give us a right action on $T_{\pi^{-1}(x)}P$
by a fixed element, but it does give a right action when acting on
the element in the identity section, so since the identity sections
remain constant we have
\begin{equation}
\begin{aligned}\vec{\varphi}_{i}^{\prime} & =\sigma_{i}^{*}\left(\gamma^{-1}\right)^{*}\vec{\varphi}_{P}\\
 & =\left(\gamma^{-1}\sigma_{i}\right)^{*}\vec{\varphi}_{P}\\
 & =\left(\gamma_{i}^{-1}\sigma_{i}\right)^{*}\vec{\varphi}_{P}\\
 & =\sigma_{i}^{*}\left(\gamma_{i}^{-1}\right)^{*}\vec{\varphi}_{P}\\
 & =\sigma_{i}^{*}\check{\gamma}_{i}\vec{\varphi}_{P}\\
 & =\check{\gamma}_{i}\vec{\varphi}_{i},
\end{aligned}
\end{equation}
where we have used $\left(g(h)\right)^{*}\varphi=h^{*}\left(g^{*}\varphi\right)$
twice, and in the penultimate line we used the equivariance of $\vec{\varphi}_{P}$.
Under neighborhood-wise gauge transformations, there is no change
in $\vec{\varphi}_{P}$ but we have new identity sections $\sigma_{i}^{\prime}(x)=\gamma_{i}^{-1}(\sigma_{i}(x))$,
so that we get
\begin{equation}
\begin{aligned}\vec{\varphi}_{i}^{\prime} & =\sigma_{i}^{\prime*}\vec{\varphi}_{P}\\
 & =\left(\gamma_{i}^{-1}\left(\sigma_{i}\right)\right)^{*}\vec{\varphi}_{P}\\
 & =\sigma_{i}^{*}\left(\gamma_{i}^{-1}\right)^{*}\vec{\varphi}_{P}\\
 & =\check{\gamma}_{i}\vec{\varphi}_{i},
\end{aligned}
\end{equation}
matching the behavior for both automorphism gauge transformations
and for gauge transformations as previously defined directly on $M$
in Section \ref{subsec:Matter-fields-and-gauges}. 

Note that if a horizontal equivariant form takes values in the abstract
fiber $F$ of another bundle associated to the frame bundle, the same
reasoning applies, but with $\check{\gamma}_{i}$ applied using the
left action of $G$ on $F$. In particular, recalling from Section
\ref{subsec:Associated-bundles} that the adjoint rep $\rho=\mathrm{Ad}$
of $G$ on $\mathfrak{g}$ defines an associated bundle $(\mathrm{Ad}P,M,\mathfrak{g})$
to $P$, we can consider a $\mathfrak{g}$-valued horizontal equivariant
form $\check{\Theta}_{P}$ on $P$, whose pullback by the identity
section under a gauge transformation satisfies 
\begin{equation}
\check{\Theta}_{i}^{\prime}=\check{\gamma}_{i}\check{\Theta}_{i}\check{\gamma}_{i}^{-1},
\end{equation}
and which similarly across trivializing neighborhoods also undergoes
a gauge transformation
\begin{equation}
\check{\Theta}_{i}=\check{g}_{ij}\check{\Theta}_{j}\check{g}_{ij}^{-1}.
\end{equation}

\section{Generalizing connections}

\subsection{Connections on bundles }

The fibers of a smooth bundle $(E,M,\pi)$ let us define vertical
tangents, but we have no structure that would allow us to canonically
define a horizontal tangent. This structure is introduced via the
\textbf{Ehresmann connection 1-form}\index{Ehresmann connection 1-form}
(AKA bundle connection 1-form\index{bundle connection 1-form}), a
vector-valued 1-form on $E$ that defines the vertical component of
its argument $v$, which we denote $v^{\baro}$, and therefore also
defines the horizontal component, which we denote $v^{\minuso}$:
\begin{equation}
\begin{aligned}\vec{\Gamma}(v) & \equiv v^{\baro},\\
H_{p} & \equiv\left\{ v\in T_{p}E\mid\vec{\Gamma}(v)=0\right\} \\
\Rightarrow v & =v^{\baro}+v^{\minuso},
\end{aligned}
\end{equation}
where $v^{\baro}\in V_{p}$, $v^{\minuso}\in H_{p}$, and $H_{p}$
is called the \textbf{horizontal tangent space}\index{horizontal tangent space}.
Viewing the $H_{p}$ as fibers over $E$ then yields the \textbf{horizontal
bundle}\index{horizontal bundle} $(HE,E,\pi_{H})$, and a \textbf{vertical
form }\index{vertical form} is defined to vanish whenever any of
its arguments are horizontal. Alternatively, one can start by defining
the horizontal tangent spaces as smooth sections of the jet bundle
of order 1 over $E$, which uniquely determines a Ehresmann connection
1-form. %
\begin{framed}%
\noindent $\triangle$ ``Ehresmann connection'' can refer to the
horizontal tangent spaces, the horizontal bundle, the connection 1-form,
or the complementary 1-form that maps to the horizontal component
of its argument.\end{framed}

Recall that on a smooth principal bundle $(P,M,\pi,G)$, the right
action $\rho\colon G\to\mathrm{Diff}(P)$ has a corresponding Lie
algebra action $\mathrm{d}\rho\colon\mathfrak{g}\to\mathrm{vect}(P)$
where $\mathrm{d}\rho\left|_{p}\right.$ is a vector space isomorphism
from $\mathfrak{g}$ to $V_{p}$. The \textbf{principal connection
1-form}\index{principal connection 1-form} (AKA principal $G$-connection\index{principal G-connection},
$G$-connection 1-form\index{G-connection 1-form}) is a $\mathfrak{g}$-valued
vertical 1-form $\check{\Gamma}_{P}$ on $P$ that defines the vertical
part of its argument $v$ at $p$ via this isomorphism, i.e. the right
action of the structure group transforms it into the Ehresmann connection
1-form: 
\begin{equation}
\begin{aligned}\mathrm{d}\rho\left(\check{\Gamma}_{P}(v)\right)\left|_{p}\right. & \equiv v^{\baro}\\
 & =\vec{\Gamma}(v)
\end{aligned}
\end{equation}
For $g\in G$, $\mathrm{d}g(v)\colon T_{p}P\to T_{g(p)}P$ preserves
horizontal tangent vectors as well as vertical. %
\begin{framed}%
\noindent $\triangle$ As with the Ehresmann connection, a ``connection''
on a principal bundle can refer to the principal connection 1-form,
the horizontal tangent spaces, or other related quantities.\end{framed}

\subsection{Parallel transport on the frame bundle}

On a frame bundle $(P=F(E),M,\pi,GL(n,\mathbb{K}))$ with connection,
we consider the horizontal tangent space to define the direction of
parallel transport. More precisely, we define a \textbf{horizontal
lift}\index{horizontal lift} of a curve $C$ from $x$ to $y$ on
$M$ to be a curve $C_{P}$ that projects down to $C$ and whose tangents
are horizontal: 
\begin{equation}
\begin{aligned}\pi\left(C_{P}\right) & =C\\
\dot{C}_{P}\left|_{p}\right. & \in H_{p}
\end{aligned}
\end{equation}
There is a unique horizontal lift of $C$ that starts at any $p\in\pi^{-1}(x)$,
whose endpoint lets us define the parallel transporter\index{parallel transporter!on the frame bundle}
on $F(E)$ 
\begin{equation}
\parallel_{C}\colon\pi^{-1}(x)\to\pi^{-1}(y).
\end{equation}
The parallel transporter is a diffeomorphism between fibers, and it
commutes with the right action: 
\begin{equation}
\parallel_{C}\left(g\left(p\right)\right)=g\left(\parallel_{C}\left(p\right)\right)
\end{equation}
We can then recover the parallel transporter on $M$ by choosing a
frame (i.e. a local trivialization), using the horizontal lift that
starts at the element $\sigma_{i}=e_{i}$ in the identity section,
and recalling the relation $e_{p}=e_{i}f_{i}(p)$:
\begin{equation}
\begin{aligned}\parallel_{C}\left(e_{i}\left|_{x}\right.\right) & =e_{i}\left|_{y}\right.f_{i}\left(\parallel_{C}\left(e_{i}\left|_{x}\right.\right)\right)\\
\Rightarrow\left(\parallel_{C}\left(v\right)\right)_{i}^{\mu}\left|_{y}\right. & =f_{i}\left(\parallel_{C}\left(e_{i}\left|_{x}\right.\right)\right)^{\mu}{}_{\lambda}v_{i}^{\lambda}\left|_{x}\right.\\
\Rightarrow\parallel^{\mu}{}_{\lambda}\left(C\right) & =f_{i}\left(\parallel_{C}\left(e_{i}\left|_{x}\right.\right)\right)^{\mu}{}_{\lambda}
\end{aligned}
\end{equation}
The second line transforms vector components using the change of basis
matrix in the opposite direction.

Similarly, on the frame bundle we can recover the connection 1-form
on $v\in T_{x}M$ within a trivializing neighborhood by using the
pullback of the identity section: 
\begin{equation}
\begin{aligned}\check{\Gamma}_{i}(v) & =\sigma_{i}^{*}\check{\Gamma}_{P}(v)\\
 & =\check{\Gamma}_{P}\left(\mathrm{d}\sigma_{i}(v)\right)
\end{aligned}
\end{equation}
On $F(E)$, $\sigma_{i}=e_{i}$ is the frame used to define the components
of vectors in the internal space on $U_{i}$, and $\check{\Gamma}_{i}(v)$
then is the element of $gl(n,\mathbb{K})$ corresponding to the vertical
component of $v$ after being mapped to a tangent of the identity
section. Thus since we consider the horizontal tangent space to define
the direction of parallel transport, $\check{\Gamma}_{i}(v)$ is the
infinitesimal linear transformation that takes the parallel transported
frame to the frame in the direction $v$, the same interpretation
found for manifolds in \cite{Marsh}. %
\begin{framed}%
\noindent $\triangle$ It is important to remember that $\check{\Gamma}_{i}$
takes values that are dependent upon the local trivialization that
defines the identity section (i.e. it is frame-dependent), while the
values of $\check{\Gamma}_{P}$ are intrinsic to the frame bundle.
This reflects the fact that the connection is a choice of horizontal
correspondences between frames, and so cannot have any value intrinsic
to $E$. \end{framed}
\begin{figure}[H]
\noindent \begin{centering}
\includegraphics[width=1\columnwidth]{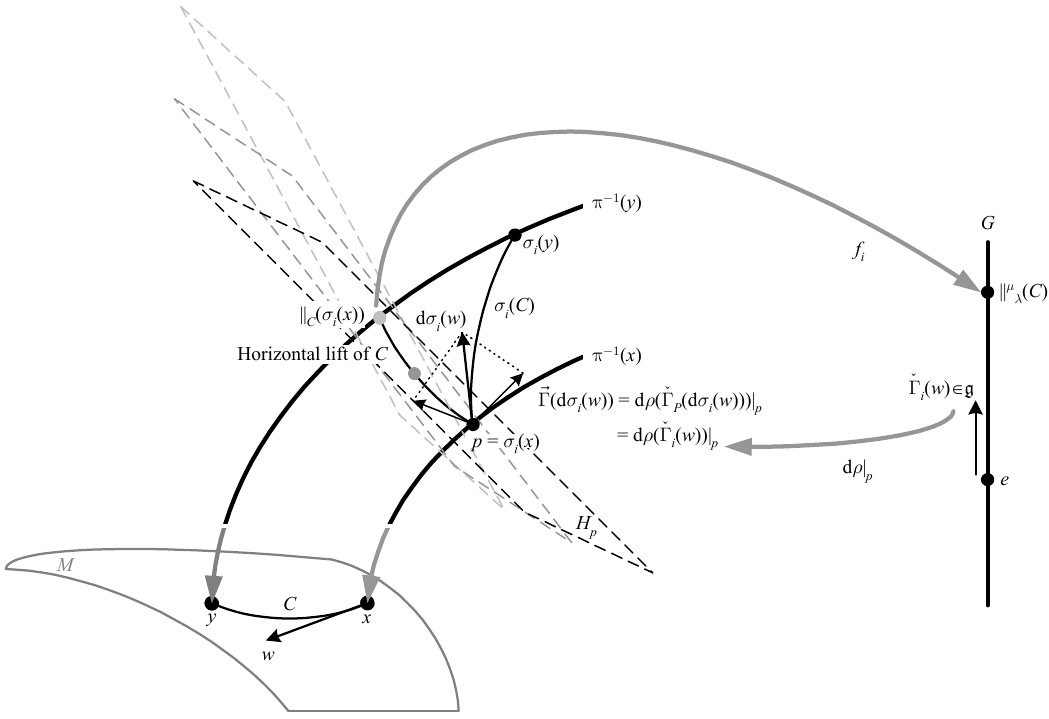}
\par\end{centering}
\caption{A principal connection 1-form on $(P,M,G)$ defines the vertical component
of its argument as a value in the Lie algebra $\mathfrak{g}$ via
the isomorphism defined by the differential of the right action $\mathrm{d}\rho$.
A horizontal lift of a curve $C$ yields the parallel transporter,
and the pullback by the identity section recovers the connection 1-form
on $M$.}
\end{figure}

The transition functions on the frame bundle can be viewed as $GL(n,\mathbb{K})$-valued
0-forms $\check{g}_{ij}$ on $U_{i}\cap U_{j}$, and it can be shown
that 
\begin{equation}
\check{\Gamma}_{i}(v)=\check{g}_{ij}\check{\Gamma}_{j}(v)\check{g}_{ij}^{-1}+\check{g}_{ij}\mathrm{d}\check{g}_{ij}^{-1}(v),
\end{equation}
which is the transformation of the connection 1-form under a change
of frame $\check{g}_{ij}^{-1}$ on manifolds. This is consistent with
the interpretation of the action of $g_{ij}$ as a change of frame
$g_{ij}^{-1}$ in Section \ref{subsec:Vector-bundles}, and it can
be shown that a unique connection on $F(E)$ is determined by locally
defined connection 1-forms on $M$ and sections that are related by
the same transition functions. %
\begin{framed}%
\noindent \sun{} The inhomogeneous transformation of the connection
1-form can be viewed as reflecting the fact that both the location
and ``shape'' of the identity section is different across local
trivializations (although we have depicted the identity sections as
``flat,'' the values of each $\sigma_{i}(x)$ are smooth but arbitrary). \end{framed}%
\begin{framed}%
\noindent \sun{} This demonstrates the advantage of the principal
bundle formulation, in that the connection 1-form on $M$ is frame-dependent,
and therefore cannot in general be defined on all of $M$, while in
contrast the principal connection 1-form is defined on all of $F(E)$,
and can be used to determine a consistent connection 1-form on $M$
within each trivializing neighborhood.\end{framed}Under either type of gauge transformation, it can also be shown that
as expected we have 
\begin{equation}
\check{\Gamma}_{i}^{\prime}(v)=\check{\gamma}_{i}\check{\Gamma}_{i}(v)\check{\gamma}_{i}^{-1}+\check{\gamma}_{i}\mathrm{d}\check{\gamma}_{i}^{-1}(v).
\end{equation}

\subsection{The exterior covariant derivative on bundles }

The exterior covariant derivative of a form on a smooth bundle with
connection is the horizontal form that results from taking the exterior
derivative on the horizontal components of all its arguments, i.e.
for a $k$-form $\varphi$ we define

\begin{equation}
\begin{aligned}\mathrm{D}\varphi(v_{0},\ldots,v_{k}) & \equiv\mathrm{d}\varphi(v_{0}^{\minuso},\ldots,v_{k}^{\minuso}).\end{aligned}
\end{equation}
On a smooth bundle, $\mathrm{D}\varphi$ can then be viewed as the
``sum of $\varphi$ on the boundary of the horizontal hypersurface
defined by its arguments.'' Note that these boundaries are all defined
by horizontal vectors except those including a Lie bracket, which
may have a vertical component. So for example, if $\varphi$ is a
vertical 1-form we have 
\begin{equation}
\mathrm{D}\varphi(v,w)=-\varphi(\left[v^{\minuso},w^{\minuso}\right]),
\end{equation}
the other terms all vanishing.

For a vector bundle $(E,M,\mathbb{K}^{n})$ associated to a smooth
principal bundle with connection $(P,M,GL(n,\mathbb{K}))$, it can
be shown that an $\mathbb{K}^{n}$-valued horizontal equivariant form
$\vec{\varphi}_{P}$ on $P$ satisfies the familiar equation

\begin{equation}
\mathrm{D}\vec{\varphi}_{P}=\mathrm{d}\vec{\varphi}_{P}+\check{\Gamma}_{P}\wedge\vec{\varphi}_{P},
\end{equation}
where the derivatives are taken on the components of $\vec{\varphi}_{P}$,
and the action of $gl(n,\mathbb{K})$-valued $\check{\Gamma}_{P}$
on the values of $\vec{\varphi}_{P}$ in is the differential of the
left action of $GL(n,\mathbb{K})$. $\mathrm{D}\vec{\varphi}_{P}$
is then also a horizontal equivariant form. Applying the pullback
by the identity section to the exterior covariant derivative, we obtain
the expected 
\begin{equation}
\mathrm{D}\vec{\varphi}_{i}=\mathrm{d}\vec{\varphi}_{i}+\check{\Gamma}_{i}\wedge\vec{\varphi}_{i}.
\end{equation}
\begin{framed}%
\noindent $\triangle$ As with the connection 1-form, it is important
to remember that the values of $\vec{\varphi}_{i}$ on $M$ are components
operated on by the matrix $\check{\Gamma}_{i}$, both of which are
defined by a local trivialization.\end{framed}

The immediate application of the above is to a $\mathbb{K}^{n}$-valued
form on the frame bundle. However, we can also apply it to other associated
bundles to $P$. In particular, recalling Section \ref{subsec:Vertical-tangents-and-horizontal-equivariant-forms},
in the associated bundle $(\mathrm{Ad}P,M,gl(n,\mathbb{K}))$ we can
apply it to a $gl(n,\mathbb{K})$-valued horizontal equivariant form
$\check{\Theta}_{P}$ on $P$, where the left action of $GL(n,\mathbb{K})$
is $\rho=\mathrm{Ad}$, and the left action of $gl(n,\mathbb{K})$
on itself is therefore $\mathrm{d}\rho=\mathrm{ad}$, i.e. the Lie
bracket. For such a form we then have 
\begin{equation}
\mathrm{D}\check{\Theta}_{P}=\mathrm{d}\check{\Theta}_{P}+\check{\Gamma}_{P}[\wedge]\check{\Theta}_{P},
\end{equation}
where again the exterior derivative is taken on the matrix components
of $\check{\Theta}_{P}$, and the action of $gl(n,\mathbb{K})$-valued
$\check{\Gamma}_{P}$ on the values of $\check{\Theta}_{P}$ is the
Lie bracket, the differential of the left action of $GL(n,\mathbb{K})$.
Applying the pullback by the identity section recovers the same formula
for algebra-valued forms on $M$ (see \cite{Marsh} for an explanation
of the notation $[\wedge]$).

\subsection{Curvature on principal bundles }

On a smooth principal bundle with connection $(P,M,G)$, the exterior
covariant derivative gives us a definition for the curvature of the
principal connection, the horizontal $\mathfrak{g}$-valued 2-form
on $P$ 
\begin{equation}
\check{R}_{P}\equiv D\check{\Gamma}_{P}.
\end{equation}
Note that the analog of the above equation on $M$ itself does not
hold. Since $\check{\Gamma}_{P}$ is vertical, this can be written
\begin{equation}
\begin{aligned}\check{R}_{P}\left(v,w\right) & =-\check{\Gamma}_{P}\left(\left[v^{\minuso},w^{\minuso}\right]\right)\\
\Rightarrow\mathrm{d}\rho\left(\check{R}_{P}\left(v,w\right)\right)\left|_{p}\right. & =-\left[v^{\minuso},w^{\minuso}\right]^{\baro},
\end{aligned}
\end{equation}
so that the curvature of the principal connection is the element of
$\mathfrak{g}$ corresponding to the vertical component of the Lie
bracket of the horizontal components of its arguments. 

On a frame bundle, we associate the horizontal tangent space with
parallel transport, and the curvature is the ``infinitesimal linear
transformation between parallel transport in opposite directions around
the boundary of the horizontal hypersurface defined by its arguments,''
or equivalently the ``infinitesimal linear transformation associated
with the vertical component of the negative Lie bracket of the horizontal
components of its arguments.'' The curvature on $M$ can be recovered
using identity sections $\sigma_{i}$ as with the connection: 
\begin{equation}
\check{R}_{i}\equiv\sigma_{i}^{*}\check{R}_{P}
\end{equation}
\begin{figure}[H]
\noindent \begin{centering}
\includegraphics[width=1\columnwidth]{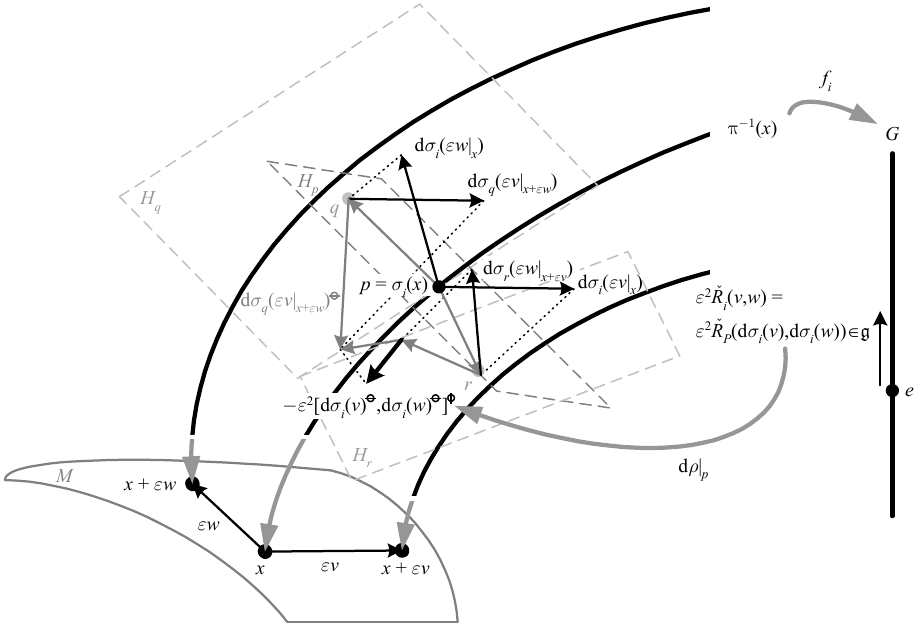}
\par\end{centering}
\caption{The curvature of the principal connection is the element of $\mathfrak{g}$
corresponding to the vertical component of the negative Lie bracket
of the horizontal components of its arguments. The sections used at
$q$ and $r$ are arbitrary, since they don't affect the vertical
component of the loop remainder. If the arguments are pulled back
using the identity section, we recover the curvature on the base space
$M$.}
\end{figure}

When $G$ is a matrix group, we find analogs of equations for curvature
on $M$ using the relations from the previous section: 
\begin{equation}
\begin{aligned}\check{R}_{P} & =\mathrm{d}\check{\Gamma}_{P}+\frac{1}{2}\check{\Gamma}_{P}[\wedge]\check{\Gamma}_{P}\\
\check{R}_{i} & =\mathrm{d}\check{\Gamma}_{i}+\frac{1}{2}\check{\Gamma}_{i}[\wedge]\check{\Gamma}_{i}\\
\mathrm{D}\check{R}_{P} & =0
\end{aligned}
\end{equation}
Note that $\check{R}_{P}$ is a map from 2-forms on $P$ to $\mathfrak{g}$,
where $G$ has a left action via the adjoint rep of $G$ on $\mathfrak{g}$.
One can then show that $\check{R}_{P}$ is equivariant with respect
to this action and that of $G$ on 2-forms, i.e. we have
\begin{equation}
g^{*}\check{R}_{P}=g_{\mathrm{Ad}}^{-1}\left(\check{R}_{P}\right).
\end{equation}
Thus $\check{R}_{P}$ is a horizontal equivariant form, and recalling
Section \ref{subsec:Vertical-tangents-and-horizontal-equivariant-forms}
we have the expected transformations 
\begin{equation}
\begin{aligned}\check{R}_{i} & =\check{g}_{ij}\check{R}_{j}\check{g}_{ij}^{-1},\\
\check{R}_{i}^{\prime} & =\check{\gamma}_{i}\check{R}_{i}\check{\gamma}_{i}^{-1}.
\end{aligned}
\end{equation}

If a flat connection (zero curvature) can be defined on a principal
bundle, then the structure group is discrete. If in addition the base
space is simply connected, then the bundle is trivial. %

\subsection{\label{subsec:The-tangent-bundle-and-solder-form}The tangent bundle
and solder form}

Returning to our motivating example, the \textbf{tangent bundle}\index{tangent bundle}
on a manifold $M^{n}$, denoted $TM$, is a smooth vector bundle $(E,M^{n},\mathbb{R}^{n})$
with a (possibly reducible) structure group $GL(n,\mathbb{R})$ that
acts as an inverse change of local frame across trivializing neighborhoods.
These trivializing neighborhoods can be obtained from an atlas on
$M$, with fiber homeomorphisms $f_{i}\colon T_{x}M\to\mathbb{R}^{n}$
defined by components in the coordinate frame $e_{i\mu}=\partial/\partial x_{i}^{\mu}$,
so that the transition functions are Jacobian matrices
\begin{equation}
\begin{aligned}v_{i}^{\mu} & =(g_{ij})^{\mu}{}_{\lambda}v_{j}^{\lambda}\\
 & =\frac{\partial x_{i}^{\mu}}{\partial x_{j}^{\lambda}}v_{j}^{\lambda}
\end{aligned}
\end{equation}
associated with the transformation of vector components. $M$ is orientable
iff these Jacobians all have positive determinant, i.e. iff the structure
group is reducible to $GL(n,\mathbb{R})^{e}$ (the definition of $TM$
being orientable). A section of the tangent bundle is a vector field
on $M$. A change of coordinates within each coordinate patch then
generates a change of frame 
\begin{equation}
\frac{\partial}{\partial x_{i}^{\prime\mu}}=\frac{\partial x_{i}^{\lambda}}{\partial x_{i}^{\prime\mu}}\frac{\partial}{\partial x_{i}^{\lambda}},
\end{equation}
which is equivalent to new local trivializations where
\begin{equation}
\left(v_{i}^{\mu}\right)^{\prime}=\frac{\partial x_{i}^{\prime\mu}}{\partial x_{i}^{\lambda}}v_{i}^{\lambda},
\end{equation}
giving us new transition functions
\begin{equation}
\frac{\partial x_{i}^{\prime\mu}}{\partial x_{j}^{\prime\lambda}}=\frac{\partial x_{i}^{\prime\mu}}{\partial x_{i}^{\sigma}}\frac{\partial x_{i}^{\sigma}}{\partial x_{j}^{\nu}}\frac{\partial x_{j}^{\nu}}{\partial x_{j}^{\prime\lambda}}.
\end{equation}

The \textbf{tangent frame bundle}\index{tangent frame bundle} (AKA
frame bundle), denoted $FM$, is the smooth frame bundle of $TM$,
i.e. $(FM,M^{n},GL(n,\mathbb{R}))$, where the fixed bases in each
trivializing neighborhood are again obtained from the atlas on $M$,
giving the same transition functions as in the tangent bundle. The
bases in $\pi^{-1}(x)$ are thus defined by 
\begin{equation}
e_{p\mu}=f_{i}(p)^{\lambda}{}_{\mu}\frac{\partial}{\partial x_{i}^{\lambda}}.
\end{equation}
A section of the frame bundle is a frame on $M$, and a global section
is a global frame, so that $M$ is parallelizable iff $FM$ is trivial.
The right action of a matrix $g^{\mu}{}_{\lambda}\in GL(n,\mathbb{R})$
operates on bases as row vectors, and an automorphism of $FM$ along
with a redefinition of fixed bases to preserve identity sections generates
changes of frame in each trivializing neighborhood that preserve the
transition functions. %
\begin{framed}%
\noindent $\triangle$ The tangent frame bundle is also denoted $F(M)$,
but rarely $F(TM)$, which is what would be consistent with general
frame bundle notation.\end{framed}

The tangent frame bundle is special in that we can relate its tangent
vectors to the elements of the bundle as bases. Specifically, we define
the \textbf{solder form}\index{solder form} (AKA soldering form,
tautological 1-form, fundamental 1-form), as a $\mathbb{R}^{n}$-valued
1-form $\vec{\theta}_{P}$ on $P=FM^{n}$ which at a point $p=e_{p}$
projects its argument $v\in T_{p}FM$ down to $M$ and then takes
the resulting vector's components in the basis $e_{p}$, i.e. 
\begin{equation}
\vec{\theta}_{P}(v)\equiv\mathrm{d}\pi(v)_{p}^{\mu}.
\end{equation}
The projection makes the solder form horizontal, and it is also not
hard to show it is equivariant, since both actions essentially effect
a change of basis:
\begin{equation}
g^{*}\vec{\theta}_{P}(v)=\check{g}^{-1}\vec{\theta}_{P}(v).
\end{equation}
The pullback by the identity section 
\begin{equation}
\begin{aligned}\vec{\theta}_{i} & \equiv\sigma_{i}^{*}\vec{\theta}_{P}\end{aligned}
\end{equation}
simply returns the components of the argument in the local basis,
and thus is identical to the dual frame $\vec{\beta}$ viewed as a
frame-dependent $\mathbb{R}^{n}$-valued 1-form. Thus recalling Section
\ref{subsec:Vertical-tangents-and-horizontal-equivariant-forms},
the values of $\vec{\theta}_{P}$ in the fiber over $x$ correspond
to a single point in the associated bundle $TM$, so that the union
of the pullbacks $\vec{\theta}_{i}$ can be viewed as a single $TM$-valued
1-form on $M$ 
\begin{equation}
\vec{\theta}\colon TM\to TM
\end{equation}
which identifies, or ``solders,'' the tangent vectors on $M$ to
elements in the bundle $TM$ associated to $FM$ (explaining the alternative
name ``tautological 1-form'').%
\begin{framed}%
\noindent $\triangle$ The $TM$-valued 1-form $\vec{\theta}$ is
also sometimes called the solder form, and can be generalized to bundles
$E$ with more general fibers as 
\begin{equation}
\theta_{E}(v)\colon TM\to E
\end{equation}
or 
\begin{equation}
\theta_{\sigma_{0}}(v)\colon TM\to V_{\sigma_{0}}E,
\end{equation}
where in the second case $\sigma_{0}$ is a distinguished section
(e.g. the zero section in a vector bundle). This is called a \textbf{soldering}
of $E$ to $M$; for example a Riemannian metric provides a soldering
of the cotangent bundle to $M$. In classical dynamics, if $M$ is
a configuration space then the solder form to the cotangent bundle
is called the Liouville 1-form\index{Liouville 1-form}, Poincaré
1-form\index{Poincaré 1-form}, canonical 1-form\index{canonical 1-form},
or symplectic potential\index{symplectic potential}.\end{framed}%
\begin{framed}%
\noindent $\triangle$ The solder form can also be used to identify
the tangent space with a subspace of a vector bundle over $M$ with
higher dimension than $M$. \end{framed}

\subsection{Torsion on the tangent frame bundle}

The covariant derivative of the solder form defines the torsion on
$P$
\begin{equation}
\begin{aligned}\vec{T}_{P} & \equiv\mathrm{D}\vec{\theta}_{P}\\
 & =\mathrm{d}\vec{\theta}_{P}+\check{\Gamma}_{P}\wedge\vec{\theta}_{P}.
\end{aligned}
\end{equation}
$\vec{T}_{P}$ is a horizontal equivariant form since $\vec{\theta}_{P}$
is. Examining the first few components, we have:
\begin{equation}
\begin{aligned}\vec{T}_{P}\left(v,w\right) & =\mathrm{d}\vec{\theta}_{P}\left(v^{\minuso},w^{\minuso}\right)\\
\Rightarrow\varepsilon^{2}\vec{T}_{P}\left(v,w\right) & =\vec{\theta}_{P}\left(\varepsilon w^{\minuso}\left|_{p+\varepsilon v^{\minuso}}\right.\right)-\vec{\theta}_{P}\left(\varepsilon w^{\minuso}\left|_{p}\right.\right)-\ldots\\
 & =\mathrm{d}\pi\left(\varepsilon w^{\minuso}\left|_{p+\varepsilon v^{\minuso}}\right.\right)_{_{p+\varepsilon v^{\minuso}}}^{\mu}-\mathrm{d}\pi\left(\varepsilon w^{\minuso}\left|_{p}\right.\right)_{p}^{\mu}-\ldots
\end{aligned}
\end{equation}
The first term projects the horizontal component of $w$ at $p+\varepsilon v^{\minuso}$
down to $M$, which is the same as projecting $w$ itself down to
$M$ since the projection of the vertical part vanishes. Then we take
its components in the basis at $p+\varepsilon v^{\minuso}$, which
is the parallel transport of the basis at $p$ in the direction $v$.
These are the same components as that of the projection of $w$ at
$p+\varepsilon v^{\minuso}$ parallel transported back to $p$ in
the basis at $p$. Thus the torsion on $P$ is the ``sum of the boundary
vectors of the surface defined by the projection of its arguments
down to $M$ after being parallel transported back to $p$.'' 

This analysis makes it clear that the pullback of the torsion on $P$
by the identity section
\begin{equation}
\vec{T}_{i}\equiv\sigma_{i}^{*}\vec{T}_{P},
\end{equation}
which by our previous pullback results recovers the torsion on $M$,
just bounces the argument vectors to the identity section and back. 

It can also be shown that the analog of the first Bianchi identity
on $M$ holds on $P$, with the original being recovered upon pulling
back by the identity section:
\begin{equation}
\begin{aligned}\mathrm{D}\vec{T}_{P} & =\check{R}_{P}\wedge\vec{\theta}_{P}\\
\mathrm{D}\vec{T}_{i} & =\check{R}_{i}\wedge\vec{\theta}_{i}
\end{aligned}
\end{equation}

\subsection{\label{subsec:Spinor-bundles}Spinor bundles}

A \textbf{spin structure}\index{spin structure} on an orientable
Riemannian manifold $M$ is a principal bundle map 
\begin{equation}
\Phi_{P}\colon(P,M^{n},\mathrm{Spin}(n))\to(OM,M^{n},SO(n))
\end{equation}
from the \textbf{spin frame bundle}\index{spin frame bundle} (AKA
bundle of spin frames\index{bundle of spin frames}) $P$ to the orthonormal
frame bundle $OM$ with respect to the double covering map 
\begin{equation}
\Phi_{G}\colon\mathrm{Spin}(n)\to SO(n).
\end{equation}
The equivariance condition on the bundle map is then 
\begin{equation}
\Phi_{P}(U(p))=\Phi_{G}(U)(\Phi_{P}(p)),
\end{equation}
so that the right action of a spinor transformation $U\in\mathrm{Spin}(n)$
on a spin basis corresponds to the right action of a rotation $\Phi_{G}(U)$
on the corresponding orthonormal basis $\Phi_{P}(p)$. %
On a time and space orientable pseudo-Riemannian manifold, a spin
structure is a principal bundle map with respect to the double covering
map $\Phi_{G}\colon\mathrm{Spin}(r,s)^{e}\to SO(r,s)^{e}$ (except
in the case $r=s=1$, which is not a double cover). 
\begin{figure}[H]
\noindent \begin{centering}
\includegraphics[width=0.9\columnwidth]{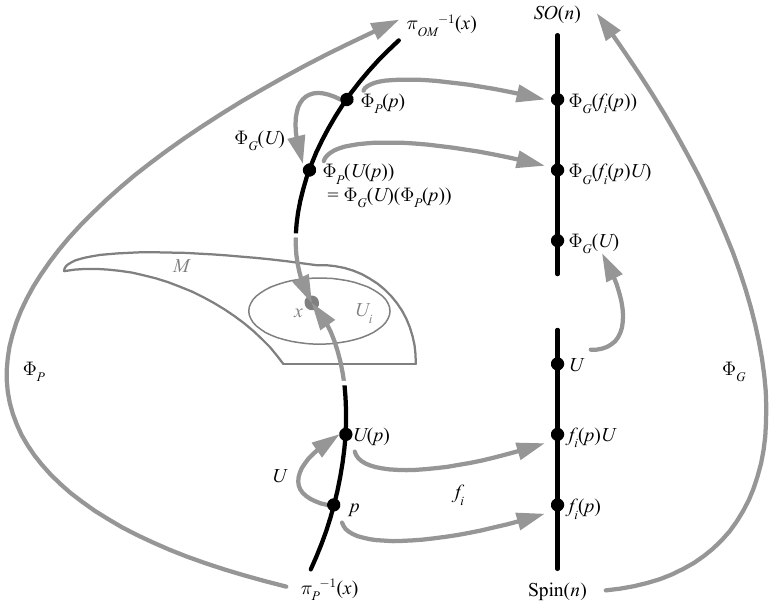}
\par\end{centering}
\caption{A spin structure is a principal bundle map that gives a global 2-1
mapping from the fibers of the spin frame bundle to the fibers of
the orthonormal frame bundle. The existence of a spin structure means
that a change of frame can be smoothly and consistently mapped to
changes of spin frame, permitting the existence of spinor fields. }
\end{figure}

If a spin structure exists for $M$, then $M$ is called a \textbf{spin
manifold}\index{spin manifold} (one also says $M$ is spin; sometimes
a spin manifold is defined to include a specific spin structure).
Any manifold that can be defined with no more than two coordinate
charts is then spin, and therefore any parallelizable manifold and
any $n$-sphere is spin. As we will see in Section \ref{subsec:Characteristic-classes},
the existence of spin structures can be related to characteristic
classes. It also can be shown that any non-compact spacetime manifold
with signature $(3,1)$ is spin iff it is parallelizable. Finally,
a vector bundle $(E,M^{n},\mathbb{C}^{m})$ associated to the spin
frame bundle $(P,M,\mathrm{Spin}(r,s)^{e})$ under a rep of $\mathrm{Spin}(r,s)^{e}$
on $\mathbb{C}^{m}$ is called a \textbf{spinor bundle}\index{spinor bundle},
and a section of this bundle is a spinor field\index{spinor field}
on $M$. 

For a charged spinor field taking values in $U(1)\otimes\mathbb{C}^{m}$,
where $\mathbb{C}^{m}$ is acted on by a rep of $\mathrm{Spin}(r,s)^{e}$,
the action of $(e^{i\theta},U)\in U(1)\times\mathrm{Spin}(r,s)^{e}$
and $(-e^{i\theta},-U)$ are identical, so that the structure group
is reducible to 
\begin{equation}
\begin{aligned}\mathrm{Spin^{\mathbb{\mathit{c}}}}(r,s)^{e} & \equiv U(1)\times_{\mathbb{Z}_{2}}\mathrm{Spin}(r,s)^{e}\\
 & \equiv\left(U(1)\times\mathrm{Spin}(r,s)^{e}\right)/\mathbb{Z}_{2},
\end{aligned}
\end{equation}
where the quotient space collapses all points in the product space
which are related by changing the sign of both components. The superscript
refers to the circle $U(1)$. A \textbf{spin}\textsuperscript{\textit{c}}\textbf{
structure}\index{spinc structure@spin\textsuperscript{\textit{c}} structure}
on an orientable pseudo-Riemannian manifold $M$ is then a principal
bundle map 
\begin{equation}
\Phi_{P}\colon(P,M^{n},\mathrm{Spin^{\mathbb{\mathit{c}}}}(r,s)^{e})\to(F{}_{SO},M^{n},SO(r,s)^{e})
\end{equation}
with respect to the double covering map 
\begin{equation}
\Phi_{G}\colon\mathrm{Spin^{\mathbb{\mathit{c}}}}(r,s)^{e}\to SO(r,s)^{e}
\end{equation}
For spinor matter fields that take values in $V\otimes\mathbb{C}^{m}$
for some internal space $V$ with structure (gauge) group $G$ with
$\mathbb{Z}_{2}$ in its center (e.g. a matrix group where the negative
of every element remains in the group), we can analogously define
a \textbf{spin}\textsuperscript{\textit{G}}\textbf{ structure}\index{spinG structure@spin\textsuperscript{\textit{G}} structure}.
It can be shown (see \cite{AvisIsham}) that spin\textsuperscript{\textit{G}}
structures exist on any four dimensional $M$ if such a $G$ is a
compact simple simply connected Lie group, e.g. $SU(2i)$; therefore
the spacetime manifold has no constraints due to spin structure in
the standard model, or in any extension that includes $SU(2)$ gauged
spinors.

\section{Characterizing bundles}

\subsection{Universal bundles}

Given a fiber bundle $(E,M,\pi,F)$ and a continuous map to the base
space 
\begin{equation}
f\colon N\to M,
\end{equation}
the \textbf{pullback bundle}\index{pullback bundle} (AKA induced
bundle\index{induced bundle}, pullback of $E$ by $f$) is defined
as 
\begin{equation}
f^{*}(E)\equiv\left\{ \left(n,p\right)\in N\times E\mid f(n)=\pi(p)\right\} ,
\end{equation}
and is a fiber bundle $(f^{*}(E),N,\pi_{f},F)$ with the same fiber
but base space $N$. Projection of $q=(n,p)\in f^{*}(E)$ onto $n$
is just the bundle projection $\pi_{f}\colon f^{*}(E)\to N$, while
projection onto $p$ defines a bundle map 
\begin{equation}
\Phi\colon f^{*}(E)\to E
\end{equation}
such that 
\begin{equation}
\pi\left(\Phi(q)\right)=f\left(\pi_{f}(q)\right)=x\in M.
\end{equation}

For any topological group $G$, there exists a \textbf{universal principal
bundle}\index{universal principal bundle} (AKA universal bundle\index{universal bundle})
$(EG,BG,G)$ such that every principal $G$-bundle $(P,M,G)$ (with
$M$ at least a CW-complex) is the pullback of $EG$ by some $f\colon M\to BG$.
The base space $BG$ is called the \textbf{classifying space}\index{classifying space}
for $G$. The pullbacks of a principal bundle by two homotopic maps
are isomorphic, and thus for a given $M$ the homotopy classes of
the maps $f$ are one-to-one with the isomorphism classes of principal
$G$-bundles over $M$. 

Every vector bundle $(E,M,\mathbb{K}^{n})$ is therefore the pullback
of the \textbf{universal vector bundle}\index{universal vector bundle}
$E_{n}(\mathbb{K}^{\infty})$ (AKA tautological bundle\index{tautological bundle},
universal bundle), the vector bundle associated to the universal principal
bundle for its structure group. It can be shown that any vector bundle
admits an inner product, so we need only consider the structure groups
$O(n)$ and $U(n)$, whose classifying spaces are each a \textbf{Grassmann
manifold}\index{Grassmann manifold} (AKA Grassmannian\index{Grassmannian})
$G_{n}(\mathbb{K}^{\infty})$. This is a limit of the finite-dimensional
Grassmann manifold $G_{n}(\mathbb{K}^{k})$, which is all $n$-planes
in $\mathbb{K}^{k}$ through the origin. Each point $x\in G_{n}(\mathbb{K}^{k})$
thus corresponds to a copy of $\mathbb{K}^{n}$, as does the fiber
over $x$ in the universal vector bundle, explaining the alternate
name ``tautological bundle.'' The total space of the associated
universal principal bundle is the \textbf{Stiefel manifold}\index{Stiefel manifold}
$V_{n}(\mathbb{K}^{\infty})$, a limit of the finite-dimensional $V_{n}(\mathbb{K}^{k})$,
defined as all ordered orthonormal $n$-tuples in $\mathbb{K}^{k}$;
the bundle projection simply sends each $n$-tuple to the $n$-plane
containing it.
\begin{figure}[H]
\noindent \begin{centering}
\includegraphics[width=0.8\columnwidth]{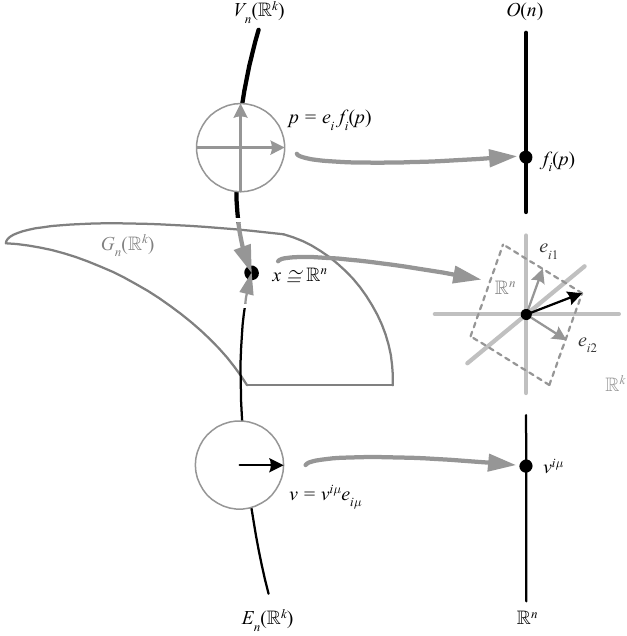}
\par\end{centering}
\caption{The Grassmann manifold $G_{n}(\mathbb{R}^{k})$ is all $n$-planes
in $\mathbb{R}^{k}$ through the origin, and is the base space of
the Stiefel manifold $V_{n}(\mathbb{R}^{k})$, defined as all ordered
orthonormal $n$-tuples in $\mathbb{R}^{k}$ where the fiber is $O(n)$
and the bundle projection simply sends each $n$-tuple to the $n$-plane
containing it. The tautological bundle is the associated vector bundle
$E_{n}(\mathbb{R}^{k})$ with fiber $\mathbb{R}^{n}$, and the universal
principal bundle for $O(n)$ is the limit $V_{n}(\mathbb{K}^{\infty})$. }
\end{figure}
\begin{framed}%
\noindent $\triangle$ Grassmann manifolds can also be denoted $Gr(n,\mathbb{K}^{k})$,
$Gr(n,k)$, $G_{n,k}$ or $g_{n,k}$ and the order of the parameters
are sometimes reversed. Stiefel manifolds have similar alternative
notations. \end{framed}

\subsection{\label{subsec:Characteristic-classes}Characteristic classes}

Vector bundles (and thus their associated principal bundles) can be
examined using \textbf{characteristic classes}\index{characteristic classes}.
For a given vector bundle $(E,M,\mathbb{K}^{n})$ these are elements
in the cohomology groups of the base space 
\begin{equation}
c(E)\in H^{\ast}(M;R),
\end{equation}
for some commutative unital ring $R$, which commute with the pullback
of any $f\colon N\to M$:
\begin{equation}
c\left(f^{*}\left(E\right)\right)=f^{*}\left(c\left(E\right)\right)
\end{equation}
In the second term, the pullback by $f$ means that $f^{*}\left(c\left(E\right)\right)\in H^{*}(N;R)$.
Since a trivial vector bundle $M\times\mathbb{K}^{n}$ is the pullback
of $(E,0,\mathbb{K}^{n})$ by $f\colon M\to0$ (where $0$ is the
space with a single point), we have 
\begin{equation}
c\left(M\times\mathbb{K}^{n}\right)=c\left(f^{*}\left(E\right)\right)=f^{*}\left(c\left(E\right)\right)=0,
\end{equation}
(where $0$ is the ring zero). Therefore the characteristic classes
of a trivial bundle vanish, or in other words a characteristic class
acts as an \textbf{obstruction}\index{obstruction} to a bundle being
trivial. However there exist non-trivial bundles whose characteristic
classes also all vanish. Similarly, if two vector bundles with the
same base space are isomorphic, then they are related by the identity
pullback; thus a necessary (but not sufficient) condition for isomorphism
is identical characteristic classes. All characteristic classes can
be determined via the cohomology classes of the classifying spaces
$BO(n)$ and $BU(n)$, since e.g. for real vector bundles any $(E,M,\mathbb{R}^{n})$
is the pullback of $BO(n)$ by some $f$, so that we have 
\begin{equation}
c\left(E\right)=c\left(f^{*}\left(BO(n)\right)\right)=f^{*}\left(c\left(BO(n)\right)\right).
\end{equation}

For a real vector bundle $(E,M,\mathbb{R}^{n})$ there are three characteristic
classes (none of which we will define here): the \textbf{Stiefel-Whitney
classes}\index{Stiefel-Whitney classes} 
\begin{equation}
w_{i}(E)\in H^{i}(M;\mathbb{Z}_{2}),
\end{equation}
the \textbf{Pontryagin classes}\index{Pontryagin classes} 
\begin{equation}
p_{i}(E)\in H^{4i}(M;\mathbb{Z}),
\end{equation}
and if the bundle is oriented the \textbf{Euler class}\index{Euler class}
\begin{equation}
e(E)\in H^{n}(M;\mathbb{Z}).
\end{equation}
For complex vector bundles, there are the \textbf{Chern classes}\index{Chern classes}
\begin{equation}
c_{i}(E)\in H^{2i}(M;\mathbb{Z}).
\end{equation}
The characteristic class of a manifold $M$ is defined to be that
of its tangent bundle, e.g. 
\begin{equation}
w_{i}(M)\equiv w_{i}(TM).
\end{equation}
If $M$ is a compact orientable four-dimensional manifold, then it
is parallelizable iff $w_{2}(M)=p_{1}(M)=e(M)=0$.

A non-zero Stiefel-Whitney class $w_{i}(E)$ acts as an obstruction
to the existence of $(n-i+1)$ everywhere linearly independent sections
of $E$. Therefore, if such section do exist, then $w_{j}(E)$ vanishes
for $j\geq i$; in particular, a non-zero $w_{n}(E)$ means there
are no non-vanishing global sections. It can be shown that $w_{1}(E)=0$
iff $E$ is orientable, so that $M$ is orientable iff $w_{1}(M)=0$. 

Spin structures exist on an oriented $M$ iff $w_{2}(M)=0$; if spin
structures do exist, then their equivalency classes have a one-to-one
correspondence with the elements of $H^{1}(M,\mathbb{Z}_{2})$. Inequivalent
spin structures have either inequivalent spin frame bundles or inequivalent
bundle maps; in four dimensions, there is only one spin frame bundle
up to isomorphism, so that different spin structures correspond to
different bundle maps (i.e. different spin connections).

Spin\textsuperscript{\textit{c}} structures exist on an oriented
$M$ if spin structures exist, but also in some cases where they do
not; for example if $M$ is simply connected and compact. If spin\textsuperscript{\textit{c}}
structures do exist, then their equivalency classes have a one-to-one
correspondence with the elements of $H^{2}(M,\mathbb{Z})$, and in
four dimensions, unlike the case for spin structures, inequivalent
spin\textsuperscript{\textit{c}} structures can have inequivalent
spin frame bundles.

\subsection{Related constructions and facts}

The tensor product\index{tensor product!vector bundles} of two vector
bundles with the same base space $(E,M,\mathbb{K}^{m})$ and $(E^{\prime},M,\mathbb{K}^{n})$
is another vector bundle 
\begin{equation}
(E\otimes E^{\prime},M,\mathbb{K}^{mn}).
\end{equation}
The direct product\index{direct product!vector bundles} of two vector
bundles $(E,M,\mathbb{K}^{m})$ and $(E^{\prime},M^{\prime},\mathbb{K}^{n})$
is another vector bundle 
\begin{equation}
(E\times E^{\prime},M\times M^{\prime},\mathbb{K}^{m+n}).
\end{equation}
If we form the direct product of two vector bundles with the same
base space, we can then restrict the base space to the diagonal via
the pullback by 
\begin{equation}
f\colon M\times M\to M
\end{equation}
defined by 
\begin{equation}
(x,x)\mapsto x.
\end{equation}
The resulting vector bundle is called the \textbf{Whitney sum}\index{Whitney sum}
(AKA direct sum bundle\index{direct sum bundle}), and is denoted
\begin{equation}
(E\oplus E^{\prime},M,\mathbb{K}^{m+n}).
\end{equation}

The \textbf{total Whitney clas}s\index{total Whitney class} of a
real vector bundle $(E,M,\mathbb{R}^{n})$ is defined as

\begin{equation}
w(E)\equiv1+w_{1}(E)+w_{2}(E)+\cdots+w_{n}(E).
\end{equation}
The series is finite since $w_{i}(E)$ vanishes for $i>n$, and is
thus an element of $H^{*}(M,\mathbb{Z}_{2})$. The total Whitney class
is multiplicative over the Whitney sum, i.e.
\begin{equation}
w(E\oplus E^{\prime})=w(E)w(E^{\prime}).
\end{equation}
The \textbf{total Chern class}\index{total Chern class} is defined
similarly, and has the same multiplicative property.

The \textbf{flag manifold}\index{flag manifold} $F_{n}(\mathbb{K}^{\infty})$
is a limit of the finite-dimensional flag manifold $F_{n}(\mathbb{K}^{k})$,
which is all ordered $n$-tuples of orthogonal lines in $\mathbb{K}^{k}$
through the origin. The name is due to the fact that an ordered $n$-tuple
of orthogonal lines in $\mathbb{K}^{k}$ is equivalent to an \textbf{\textit{n}}\textbf{-flag}\index{n-flag},
a sequence of subspaces $V_{1}\subset\cdots\subset V_{n}$ in $\mathbb{K}^{k}$
where each $V_{i}$ has dimension $i$.


\addcontentsline{toc}{section}{References}
\begin{thebibliography}{1}
\bibitem{AvisIsham} S. J. Avis and C. J. Isham, ``Generalized Spin
Structures on Four Dimensional Space-Times,'' Commun. Math. Phys.
72 (1980), 103-118, https://projecteuclid.org/euclid.cmp/1103907653

\bibitem{FrankelGeom} T. Frankel, \textit{The Geometry of Physics}
(Cambridge University Press, 1997)

\bibitem{GockelerSchucker} M. Göckeler and T. Schücker, \textit{Differential
Geometry, Gauge Theories, and Gravity} (Cambridge University Press,
1987)

\bibitem{Kobayashi}S. Kobayashi and K. Nomizu, \textit{Foundations
of Differential Geometry} (John Wiley \& Sons, 1963)

\bibitem{Lawson}H. B. Lawson, Jr. and M. Michelsohn, \textit{Spin
Geometry} (Princeton University Press, 1989)

\bibitem{Marsh}A. Marsh, ``Riemannian Geometry: Definitions, Pictures,
and Results,'' arXiv:1412.2393 {[}gr-qc{]}, http://arxiv.org/abs/1412.2393
\end{thebibliography}
\end{document}